\documentclass[12pt]{article} 
\usepackage{amsmath,amssymb}
\usepackage[dvipdfmx]{graphicx}
\usepackage{mathdots}   

\allowdisplaybreaks[1]   

\newcommand{\Section}[1]{%
\renewcommand{\thesection}{\S\arabic{section}}
\section{#1}
\renewcommand{\thesection}{\arabic{section}}
\setcounter{equation}{0}}
\newcommand{\qed}{\hfill \hbox{\rule[-2pt]{4pt}{7pt}}}   
\newcommand{\proof}{{\hspace*{0.4cm} {\it Proof}.\ \enskip}}  

\newtheorem{lem}{Lemma}[section]
\newtheorem{thm}[lem]{Theorem}  
\newtheorem{prop}[lem]{Proposition}

\newtheorem{rem}[lem]{Remark}

\newcommand{\Q}{{\Bbb Q}}
\newcommand{\Z}{{\Bbb Z}}

\newcommand{\N}{{\Bbb N}}

\newcommand{\calF}{{\cal F}}
\newcommand{\calG}{{\cal G}}

\newcommand{\calO}{{\cal O}}
\newcommand{\calP}{{\cal P}}

\newcommand{\calS}{{\cal S}}

\newcommand{\pair}[2]{\left\langle {#1}, \, {#2}\right\rangle}

\newcommand{\set}[2]{\left\{\left.#1\vphantom{#2}\:\right\vert\:#2\right\}}
\newcommand{\wt}{\widetilde}
\newcommand{\what}{\widehat}

\newcommand{\frp}{{\frak p}}

\newcommand{\lam}{{\lambda}}
\newcommand{\Lam}{{\Lambda}}

\newcommand{\alp}{{\alpha}}  
\newcommand{\vphi}{{\varphi}} 
\newcommand{\eps}{{\epsilon}} 
\newcommand{\abs}[1]{\left\vert{#1}\right\vert}  

\newcommand{\slit}{\vspace{5mm}\noindent}
\newcommand{\mslit}{\vspace{3mm}\noindent}
\newcommand{\mmslit}{\vspace{1mm}\noindent}   
\newcommand{{\any}}{{}^\forall}
\newcommand{{\is}}{{}^\exists}
\newcommand{{\st}}{\; {\rm s.t.}\; }
\newcommand{\ol}[1]{\overline{#1}}

\newcommand{\twomatrix}[4]{\begin{pmatrix}
                           {#1} & {#2}\\
                           {#3} & {#4}
                          \end{pmatrix}}
\newcommand{\twomatrixplus}[4]{\left(\begin{array}{c|c}
                           {#1}  & {#2}\\
                           \hline
                           {#3} &  {#4}
                          \end{array}\right)}

\newcommand{\calo}{{\mathfrak o}}   
\newcommand{\Trd}{{\rm T_{rd}}}    
\newcommand{\Nrd}{{\rm N_{rd}}}   
\newcommand{\gen}[1]{\langle{#1}\rangle}

\title{Local densities of $p$-adic quaternion hermitian forms}
\author{{\Large Yumiko Hironaka}\\[2mm]
Waseda University}
\begin{document}
\date{}
\maketitle

\vspace{1cm}
\begin{center}
\begin{minipage}{15cm}
\begin{center}{\bf abstract} \end{center}

We consider local densities for $p$-adic quaternion hermitian forms (hermitian forms over a division quaternion algebra over a ${\mathfrak p}$-adic field $k$). The author has studied such forms in connection with spherical functions on the space of quaternion hermitian forms in the previous paper.
Obtaining good explicit formulas of local densities is an interesting but difficult problem in general, and without explicit formulas we may study some general theory on local densities. In this paper, we give two results on the relations among local densities.
First, regarding the local density $\mu(B, A)$ (of representations of $B$ by $A$) as functions of $B$, we study their linear independence when $A$ varies by scaling hyperbolic planes.
Then, following Kitaoka for symmetric forms, we introduce a formal power series $P(B, A; X) = \sum_{r \geq 0}\, \mu(\pi^rB, A)X^r$, where $\pi$ is a prime element of $k$. We give an explicit polynomial $f(X)$ for which $f(X)\cdot P(B, A; X)$ becomes a polynomial of $X$.
Both results are based on Gauss sum expressions of local densities and Fourier transforms of functions on the space of quaternion hermitian forms.  Finally, as an appendix, we list the previous results on symmetric forms and hermitian forms for comparison.

\slit


\end{minipage}
\end{center}

\vspace{1cm}
\setcounter{section}{-1}

\setcounter{section}{-1}
\Section{Introduction}
Spherical functions on a $p$-adic homogeneous space are one of the basic objects to study the space. The cases of sesquilinear forms are particular interesting, since spherical functions can be regarded as generating functions of local densities, and the latter have a close connection to the global theory of automorphic forms.

\mmslit
In the previous paper \cite{Q}, we have introduced the space $X$ of $p$-adic quaternion hermitian forms, and studied harmonic analysis on $X$ by using spherical functions on $X$ as a main tool. 
In the cases of alternating forms and unramified hermitian forms, we gave the explicit formulas of spherical functions whose main part consists of specialized Hall-Littlewood symmetric polynomials. Since they are well studied and many relations among them are known, we could extract local densities of forms from spherical functions (cf. \cite{JMSJ}, \cite{HS-Am}). 
For the case of quaternion hermitian forms, we gave the explicit formula of spherical functions  by using certain symmetric polynomials. Unfortunately, we do not know well about these symmetric polynomials, so we have not succeeded to have good expressions of local densities. 

\mmslit
One has most interest in the case of symmetric forms, since local densities appear in Fourier coefficients of genus theta series and have a close relation to the study of Siegel modular forms (e.g., \cite{Ik}, \cite{Ka}).  Unfortunately, neither good explicit formulas of spherical functions nor local densities are known. But without explicit formulas, one may study certain general theory of local densities, and relate it to the global study of modular forms.  
The linear independence of certain local densities is closely related to the basis problem of Siegel modular forms (cf. \cite{BS}, \cite{BHS}, \cite{Kitaoka}, \cite{Kitaoka2}, \cite{Ki}, \cite{K-SP}).

\mslit
In this paper, we study similar problems for the case of quaternion hermitian forms. 
\cite{Q} is the basic reference for quaternion hermitian forms and we follow the notation there. 
In the following we introduce our main results, for the detailed definitions and notations, see \S1 and the beginning of \S2.       
Let $D$ be a division quaternion algebra over a $\frp$-adic field $k$ of odd residual characteristic, and denote by $\calO$ its ring of integers. 
By using the standard $k$-involution $*$ on $D$, we define the conjugate transpose $A^* \in M_{nm}(D)$ of $A = (a_{ij}) \in M_{mn}(D)$ by ${A^*}_{ij} = {a_{ji}}^* \in D$ for any $i, j$. 
Then the spaces $X_n$ and $X_n^+$ of quaternion hermitian forms are given as follows.
\begin{eqnarray*}
X_n = \set{x = (x_{ij}) \in GL_n(D)}{x= x^*} \supset X_n^+ = X_n \cap M_n(\calO).
\end{eqnarray*}
For $A \in X_m^+$ and $B \in  X_n^+$ with $m \geq n$,   we define the local density $\mu(B, A)$ of $B$ by $A$ as follows.
\begin{eqnarray*}
\mu(B, A) &=& \lim_{\ell \rightarrow \infty}
\dfrac{
\sharp\set{\overline{v} \in M_{mn}(\calO/\pi^\ell\calO)}{v^*Av - B \in M_n(\pi^\ell\calO)}}
{q^{\ell{n(4m-2n-1)}}}, 
\end{eqnarray*}
where $\pi$ is a prime element of $k$ and $q$ is the number of the residue class field of $k$. We may consider $\mu(B, A)$ for all $A \in X_m$ and $B \in X_n$ by scaling. 
Further we set 
\begin{eqnarray*}
&&
\Gamma_k = \set{\lam \in \Z^k}{\lam_1 \geq \cdots \geq \lam_k \geq 0} \supset \Gamma_{k, \ell} = \set{\lam \in \Gamma_k^+}{\ell \geq \lam_1}, \\
&&
\mbox{for } \lam \in \Gamma_k,  \nonumber\\
&& 
\qquad H^\lam = {h}^{\lam_1}\bot {h}^{\lam_2} \bot \cdots \bot {h}^{\lam_k} \in X_{2k}^+, \quad h^r = \twomatrix{0}{\varPi^r}{(-\varPi)^r}{0},\nonumber \\ 
 && \label{H-lam}
\qquad 
H^{{\bf 0}, k} \in X_{2k}^+, \; \mbox{for }\lam= {\bf 0} = (0, \ldots, 0) \in \Gamma_k.
\end{eqnarray*}
For $\lam \in \Gamma_{k, \ell}$, we define $\what{\lam} \in \Gamma_{\ell, k}$ by $\what{\lam}_i = \sharp\set{j}{\lam_j \geq i}, \; 1 \leq i \leq \ell$. For $\lam, \mu \in \Gamma_k$, we set $\pair{\lam}{\mu} = \sum_{i=1}^k \lam_i\mu_i$.
Then, our first main result is the following.

\slit
{\bf Theorem 1} (Theorem~2.1) {\it 
Assume that $k, \ell, n \in \N$, and $r \in \Z, r \geq 0$ satisfy the condition 

\begin{center}
$k \geq n$ \; and \; $2k+r \geq 8n-1$.
\end{center}

\noindent
Take an $S \in X_r^+$ if $r > 0$. Then one has the following.

\mmslit
{\rm (1)} For $\lam \in \Gamma_{k,\ell}$, $T \in X_n^+$, it satisfies
\begin{eqnarray*}
\mu(T, H^\lam \bot S) = 
\sum_{\tau \in \Gamma_{n, \ell}}\, a_\tau q^{2\pair{\what{\lam}}{\what{\tau}}}
\end{eqnarray*}
where, $a_\tau= a_\tau(k, S, T)$ is a rational constant independent of $\lam$.

\mmslit
{\rm (2)} As functions of $T \in X_n^+$, the set $\set{\mu(T, H^\lam \bot S)}{\lam \in \Gamma_{k,\ell}}$ spans a $\begin{pmatrix}n+\ell\\n \end{pmatrix}$-dimensional $\Q$-space, and the set $\set{\mu(T, H^\mu\bot H^{{\bf0},k-n} \bot S)}{\mu \in \Gamma_{n,\ell}}$ forms a basis.
}

\slit
Next, in \S4 we define a formal power series (Kitaoka series) for $A \in X_m$ and $B \in X_n$ with $m \geq n$ by
\begin{eqnarray*}
P(B, A; X) = \sum_{r \geq 0}\, \mu(\pi^rB, A)X^r.
\end{eqnarray*}
Then, our second main result is the following.

\slit
{\bf Theorem 2} (Theorem~4.1) {\it 
Assume that $A \in X_m$ and $B \in X_n$ with $m \geq n$. Then 
\begin{eqnarray*}
P(B, A; X) \times \prod_{i=0}^n\, (1 - q^{(n-i)(2n+2i-2m-1)}X)
\end{eqnarray*}
is a polynomial in $X$.
}

\mslit
Gauss sums and Fourier transforms are main tools to study these problems on local densities. Under some conditions, local densities can be viewed as Fourier transforms of Gauss sums (Proposition~2.4). Then  we calculate the Gauss sums explicitly (Proposition~3.2) and give the proof of Theorem\,1 in \S3. In \S4, after defining Kitaoka series, we give a modified expression of local densities by using the finite Gauss sums and finite Fourier transforms (Proposition~4.3), and introduce partial series of Kitaoka series (Proposition~4.6 and (4.18)). Then we calculate the partial series and complete the proof of Theorem\,2 in \S5.
Finally in \S6, as an appendix, we list the previous results on symmetric forms and hermitian forms for comparison.

\vspace{2cm}
\Section{Notations and preliminaries}
\mslit
Let $k$ a $\frp$-adic field, and set $\calo = \calO_k, \; \pi= \pi_k, \; \frp = \pi\calo, \; q = \abs{\calo/\frp}$. We assume that $q$ is odd. 
Let $D$ be a division quaternion algebra over $k$, and we fix related notations as follows.  
A standard $k$-basis $\{1, \eps, \varPi, \eps\varPi\}$ of $D$, such that it forms an $\calo$-basis for $\calO= \calO_D$, $k(\eps)/k$ is unramified, $\varPi^2 = \pi$, and $\eps\varPi = -\varPi\eps$.  We denote by $\calP = \varPi\calO$ the maximal ideal of $\calO$. 
We define a $k$-involution $*$ on $D$ by
\begin{eqnarray}   \label{shape alp}
\alp = a+b\eps+c\varPi+d\varPi\eps \longmapsto \alp^* = a-b\eps-c\varPi-d\varPi\eps, \; (a,b,c,d \in k),
\end{eqnarray}
where it satisfies that $(\alp\beta)^* = \beta^*\alp^*$ for $\alp, \beta \in D$.
For $A \in M_{mn}(D)$, we define its conjugate transpose $A^*$ by 
\begin{eqnarray} \label{conj-tr}
A^* \in M_{nm}(D), \quad (A^*)_{ij} = (A_{ji})^* \; (\any i, j).
\end{eqnarray}
We define the space of hermitian forms by
\begin{eqnarray}
V_n &=& V_n(D) = \set{x \in M_n(D)}{x^*=x} \supset  X_n=V_n(D)\cap GL_n(D), \nonumber\\
V_n(\calO) &=& V_n(D) \cap M_n(\calO) \supset  X_n^+ = X_n \cap M_n(\calO). \label{spaces}
\end{eqnarray}
The group $GL_n(D)$ acts on $V_n$ by 
\begin{eqnarray}
g\cdot x = x[g^*]=gxg^*, \quad g\in GL_n(D), \; x \in V_n.
\end{eqnarray}
It is clear that $X_n$ is stable under the action of $G_n=GL_n(D)$, and $X_n^+$ is stable under the action of $K_n=GL_n(\calO)$.

\mmslit
For $A \in X_m^+$ and $B \in  X_n^+$ with $m \geq n$,   we define local density of $B$ by $A$ as follows:
\begin{eqnarray}
\mu(B, A) &=& \lim_{\ell \rightarrow \infty}
\dfrac{
\sharp\set{\overline{v} \in M_{mn}(\calO/\pi^\ell\calO)}{A[v] - B \in M_n(\pi^\ell\calO)}}
{q^{\ell{n(4m-2n-1)}}}, \label{loc-dens2}
\end{eqnarray}
where we identify $M_{mn}(\calO/\pi^\ell\calO)$ with $M_{mn}(\calO)/M_{mn}(\pi^\ell\calO)$,  and $A[v] = v^*Av \in V_n(\calO)$.  If $\ell$ is big enough, the above ratio is stable (hence it is a rational number), and for that it is enough to assume that 
\begin{eqnarray} \label{assump}
\mbox{The $K_n$-orbit containing $B$ decomposes into cosets modulo $V_n(\pi^\ell\calO)$.}
\end{eqnarray}
For the calculation of local densities, it is convenient to introduce another expression
\begin{eqnarray}
&& \mu(B, A) = \lim_{\ell \rightarrow \infty}
\dfrac{N_\ell(B, A)}{q^{\ell n(4m-2n-1)-n(n-1)}}, \nonumber \\
&&
\quad
N_\ell(B, A) = \sharp \set{\overline{v} \in M_{mn}(\calO/\pi^\ell\calO)}{A[v] - B \in V_n(\pi, \ell)} \nonumber\\
&&
\quad
V_n(\pi, \ell) = \set{y \in V_n}{y_{ii} \in \frp^\ell, \; y_{ij} \in \calP^{2\ell-1}\, \;(\any i, j )}. \label{loc-dens3}
\end{eqnarray}
By the relation (\cite[Prop.2.4]{Q})
\begin{eqnarray} \label{shift}
\mu(\pi^eB, \pi^eA) = q^{en(2n-1)}\mu(B, A),
\end{eqnarray}
we may consider local density $\mu(B, A)$ for all $A \in X_m$ and $B \in X_n$.

\slit
There is a $k$-algebra map $\vphi: D \longrightarrow M_2(k(\eps))$ determined by the equation
\begin{eqnarray} \label{vphi}
\alp (1 \varPi) = (1 \varPi)\vphi(\alp) \quad (\mbox{in } D \oplus D).
\end{eqnarray}
For $\alp$ of the shape \eqref{shape alp}, we see that 
\begin{eqnarray} \label{vphi(alp)}
\vphi(\alp)  = \twomatrix{a+b\eps}{(c-d\eps)\pi}{c+d\eps}{a-b\eps} \in M_2(k(\eps)).
\end{eqnarray}
Then, the reduced norm $\Nrd$ and the reduced trace $\Trd$ on $D$ are defined by
\begin{eqnarray}
\Nrd(\alp) &=& \alp\alp^* = \det(\vphi(\alp)) \; \in k, \nonumber \\
\Trd(\alp) &=& \alp + \alp^* = {\rm trace}( \vphi(\alp)) \; \in k.
\end{eqnarray}
It is easy to see the following properties of $\Nrd$ and $\Trd$: 
\begin{eqnarray}
&&
\Nrd(\calO^\times) = \calo^\times, \quad \Trd(\calP^m) = \frp^{[\frac{m+1}{2}]} \;\; (m \in \Z), \nonumber \\
&&
\Trd(\alp\beta) = \Trd(\beta\alp). \label{Nrd-Trd}
\end{eqnarray}
We define a pairing $\pair{\;}{\;}: V_n \times V_n \longrightarrow k$ by
\begin{eqnarray} \label{pairing}
\pair{A}{B} = \sum_{i=1}^n\,  A_{ii}B_{ii} + \sum_{i < j}\, \Trd(A_{ij}B_{ji}).
\end{eqnarray}
We extend the above map $\vphi = \vphi_1$ (cf. \eqref{vphi(alp)})  to a $k$-algebra map $\vphi_n$ 
\begin{eqnarray}
\vphi_n &:& M_n(D) \longrightarrow M_{2n}(k(\eps)), \; A= (a_{ij}) \longmapsto \left(\vphi(a_{ij})\right)_{1 \leq i,j \leq n}, \label{vphi-n}
\end{eqnarray}
where we may consider the pairing by using matrices over the field $k(\eps)$. 

\begin{lem} \label{lem 1.1}
Assume that $A, B \in V_n$ and $U \in M_n(D)$. Then, one has the following.
\begin{eqnarray}
&&
{\rm (1)} \; \pair{A}{B} = \pair{B}{A} \in k. \nonumber \\
&&
{\rm (2)} \; \pair{A}{B} = \frac12{\rm trace }(\vphi_n(A)\vphi_n(B)).\nonumber \\
&&
{\rm (3)} \; \pair{A[U]}{B} = \pair{A}{B[U^*]}. \nonumber 
\end{eqnarray}
\end{lem}
\proof
(1) The fact $\pair{A}{B} \in k$ is clear by definition. By \eqref{Nrd-Trd}, we have $\Trd(\alp\beta) = \Trd(\beta^*\alp^*) = \Trd(\alp^*\beta^*)$, thus we see that
\begin{eqnarray*} 
\pair{A}{B} &=& \sum_{i=1}^n\, A_{ii}B_{ii} + \frac12\sum_{i<j}(\Trd(A_{ij}B_{ji})+\Trd(A_{ji}B_{ij})) \\
&=&
\frac12\sum_{i, j}\, \Trd(A_{ij}B_{ji}) = \pair{B}{A}.
\end{eqnarray*}

\mmslit
(2) It is clear for size $1$, since $\vphi_1(a) = \twomatrix{a}{0}{0}{a}$. 
Now we assume that the identity (2) holds for size $n$, and consider the case of size $n+1$. Set
\begin{eqnarray*}
\wt{A} = \twomatrixplus{A}{ \begin{matrix}\alp_1\\\vdots\\\alp_n\end{matrix}}
{\alp_1^* \cdots \alp_n^*}{a_0}, \quad
\wt{B} = \twomatrixplus{B}{ \begin{matrix}\beta_1\\\vdots\\\beta_n\end{matrix}}
{\beta_1^* \cdots \beta_n^*}{b_0} \in V_{n+1}(D).
\end{eqnarray*}
Then, by definition, 
\begin{eqnarray*}
&&
{\rm trace}(\vphi_{n+1}(\wt{A})\vphi_{n+1}(\wt{B})) \\
&=&
{\rm trace}(\vphi_{n}(A)\vphi_{n}(B)) + 2a_0b_0 + \sum_{i=1}^n {\rm trace}(\vphi_1(\alp_i)\vphi_1(\beta_i^*)+\vphi_1(\alp_i^*)\vphi_1(\beta_i))\\
&=&
2\left(\pair{A}{B}+a_0b_0 + \sum_{i=1}^n \Trd(\alp_i\beta_i^*)\right) \\
&=& 2\pair{\wt{A}}{\wt{B}},
\end{eqnarray*} 
which proves the identity (2) for size $n+1$.

\mmslit
(3) By the identity (2) and a property of the trace of matrices over a field, we have for $U \in M_n(D)$,
\begin{eqnarray}
2\pair{A[U]}{B} &=& {\rm trace}(\vphi_n(U^*AU)\vphi_n(B)) = {\rm trace}(\vphi_n(U^*)\vphi_n(A)\vphi_n(U)\vphi_n(B)) \nonumber\\
&=& {\rm trace}(\vphi_n(A)\vphi_n(UBU^*)) = 2\pair{A}{B[U^*]}, \nonumber
\end{eqnarray}
which completes the proof.
\qed

\mslit
\begin{align}
&
\mbox{We take and fix an additive character $\psi$ on $k$ of conductor $\calo$. For any $\ell>0$, we define the}\nonumber\\ 
& 
\mbox{character $\psi_\ell$ of conductor $\frp^\ell$ by $\psi_\ell(x) = \psi(\pi^{-\ell}x)$, and denote by $\chi_\ell$ the induced character}\nonumber\\
& \mbox{on $\calo/\frp^\ell$, which is nontrivial on $\frp^{\ell-1}/\frp^\ell$.} \label{character}
\end{align}

\mslit
Then, we have the following integral expression for $\mu(B, A)$.
\begin{lem} \label{lem 1.2} 
Let $A \in X_m$ and $B \in X_n$ with $m \geq n$. Then
\begin{eqnarray*}
\mu(B, A) = q^{-n(n-1)} \int_{V_n(D)}\, dy \int_{M_{mn}(\calO)}\, \psi(\pair{A[v]-B}{y}) dv,
\end{eqnarray*}
where the integral over $V_n(D)$ is understood as
\begin{eqnarray*}
\int_{V_n(D)} = \lim_{\ell \rightarrow \infty}\, \int_{V_n(\calP^{-2\ell})},
\end{eqnarray*}
and $dv$ and $dy$ are Haar measure on $M_{mn}(D)$ and $V_n$, respectively, normalized by
\begin{eqnarray*}
\int_{M_{mn}(\calO)} dv = 1, \quad \int_{V_n(\calO)} dy = 1.
\end{eqnarray*}
\end{lem}

\proof
We may assume that $A$ and $B$ are integral (cf. \eqref{shift}). Since $A[v]-B \pmod{V(\calP^{2\ell})}$ is stable on each coset of $M_{mn}(\calO/\calP^{2\ell})$ for $\ell >0$, we have
\begin{eqnarray}
\lefteqn{\int_{V_n(\calP^{-2\ell})}\, dy \int_{M_{mn}(\calO)}\, \psi(\pair{A[v]-B}{y})dv}\nonumber\\
&=&
q^{\ell n(2n-1)} \int_{V_n(\calO)}\, dy \int_{M_{mn}(\calO)}\, \psi(\pi^{-\ell}\pair{A[v]-B}{y})dv \nonumber\\
&=&
q^{\ell n(2n-1) -4\ell mn} \int_{V_n(\calO)} dy \sum_{\ol{v} \in M_{mn}(\calO/\calP^{2\ell})}\,  \psi_\ell(\pair{A[v]-B}{y}). \label{cal mid}
\end{eqnarray}
Set $A[v]-B = (c_{ij})$ for a while and recall the definition \eqref{pairing}. Then, by \eqref{Nrd-Trd}  and the orthogonal property of characters, we see that \eqref{cal mid} is equal to  
\begin{eqnarray}
&&
q^{-\ell n(4m-2n+1)} \sum_{\ol{v} \in M_{mn}(\calO/\calP^{2\ell})}\, \prod_{i=1}^n\, \int_{\calo} \psi_\ell(c_{ii}y)dy\, \prod_{i<j}\, \int_{\calO}\, \psi_\ell(\Trd(c_{ij}y))dy \nonumber\\
&=&
q^{-\ell n(4m-2n+1)} \sharp\set{\ol{v}\in M_{mn}(\calO/\calP^{2\ell})}{A[v]-B \in V_n(\pi, \ell)} \nonumber\\
&=&
q^{-\ell n(4m-2n+1)}N_\ell(B,A). \label{cal end}
\end{eqnarray} 
If $\ell$ is big enough, for example \eqref{assump} is satisfied, then \eqref{cal end} is equal to $q^{n(n+1)} \mu (B, A)$, which completes the proof. 
\qed

\slit
We may rearrange the integral over $V_n(\calO)$ in \eqref{cal mid}, and we obtain the following identity (cf. \cite[Prop.2.7]{Q}), where we use the symbol $\chi_\ell$ as a character on $\calo/\frp^\ell$.    
\begin{eqnarray} \label{N in finite}
N_\ell(B, A) = q^{-\ell n(2n-1)} \sum_{\ol{y} \in V_n(\calO/\calP^{2\ell})}\, \sum_{\ol{v} \in M_{mn}(\calO/\calP^{2\ell})}\, \chi_\ell(\pair{A[v]-B}{y}).
\end{eqnarray}

\slit
Before closing this section we note the $K_n$-orbit decomposition of $X_n$, which is given in \cite[Th.6.2]{Jac}. 
The set $K_n \backslash X_n$ bijectively corresponds to the set 
\begin{eqnarray}\label{Lam}
\Lam_n = \set{ \gamma \in \Z^n}{ \begin{array}{l}
\gamma_1 \geq \cdots \geq \gamma_n\\
\mbox{if $\gamma_i$ is odd, then $\sharp\set{j}{\gamma_j = \gamma_i}$ is even}\end{array}}.
\end{eqnarray}
In fact, for each $\gamma \in \Lam_n$, 
writing $\gamma$ as 
\begin{eqnarray*}
\gamma = (\lam_1^{r_1}\cdots \lam_t^{r_t}), \quad \lam_1 > \cdots > \lam_t, \; r_i > 0, \; \sum_i r_i = n,
\end{eqnarray*}
we associate an element $\pi^\gamma$ in $X_n$ as follows:
\begin{eqnarray}
&& \label{pi-gamma}
\pi^\gamma = \pi^{{\lam_1}^{r_1}} \bot \cdots \bot \pi^{{\lam_1}^{r_1}},\\
&&
\quad \pi^{\lam^r} = \left\{\begin{array}{ll}
 Diag(\pi^e, \ldots, \pi^e) & \mbox{if }\lam=2e \\
\twomatrix{0}{\pi^e\varPi}{-\pi^e\varPi}{0} \bot \cdots \bot \twomatrix{0}{\pi^e\varPi}{-\pi^e\varPi}{0} & \mbox{if }\lam=2e+1
\end{array}\right\} \in X_r. \nonumber
\end{eqnarray}
Then, the set $K_n\backslash X_n^+$ bijectively corresponds to the set
\begin{eqnarray}\label{Lam+}
\Lam_n^+ = \set{ \lam \in \Lam_n}{\lam_n \geq 0}.
\end{eqnarray}
\vspace{2cm}
\setcounter{section}{1}
\Section{Linear independence of local densities}
We consider local densities $\mu(B,A)$ as functions of $B \in X_n^+$ and study their linear independence when $A \in X_m^+$ varies by scaling hyperbolic planes. We have studied in \cite{BHS} the same problem on local densities of symmetric forms and applied it to the global study of Siegel modular forms. Here we concentrate in the local theory.

\mslit
To state our main results, we need some more notation.
\begin{eqnarray}
&&
\Gamma_k = \set{\lam \in \Z^k}{\lam_1 \geq \cdots \geq \lam_k \geq 0} \supset \Gamma_{k, \ell} = \set{\lam \in \Gamma_k^+}{\ell \geq \lam_1}, \\
&&
\mbox{for } \lam \in \Gamma_k,  \nonumber\\
&& 
\qquad H^\lam = {h}^{\lam_1}\bot {h}^{\lam_2} \bot \cdots \bot {h}^{\lam_k} \in X_{2k}^+, \quad \quad h^r = \twomatrix{0}{\varPi^r}{(-\varPi)^r}{0},\nonumber \\ 
 && \label{H-lam}
\qquad 
H^{{\bf 0}, k} \in X_{2k}^+, \; \mbox{for }\lam= {\bf 0} = (0, \ldots, 0) \in \Gamma_k.
\end{eqnarray}
For $\lam \in \Gamma_{k, \ell}$, we define $\what{\lam} \in \Gamma_{\ell, k}$ by $\what{\lam}_i = \sharp\set{j}{\lam_j \geq i}, \; 1 \leq i \leq \ell$. For $\lam, \mu \in \Gamma_k$, we set $\pair{\lam}{\mu} = \sum_{i=1}^k \lam_i\mu_i$.

\mslit
We introduce our main results.

\begin{thm} \label{th: lin indep}
Assume that $k, \ell, n \in \N$, and $r \in \Z, r \geq 0$ satisfy the condition 

\begin{center}
$k \geq n$ \; and \; $2k+r \geq 8n-1$.
\end{center}

\noindent
Take an $S \in X_r^+$ if $r > 0$. Then one has the following.

\mmslit
{\rm (1)} For $\lam \in \Gamma_{k,\ell}$, $T \in X_n^+$, it satisfies
\begin{eqnarray*}
\mu(T, H^\lam \bot S) = 
\sum_{\tau \in \Gamma_{n, \ell}}\, a_\tau q^{2\pair{\what{\lam}}{\what{\tau}}},
\end{eqnarray*}
where, $a_\tau= a_\tau(k, S, T)$ is a constant independent of $\lam$.

\mmslit
{\rm (2)} As functions of $T \in X_n^+$, the set $\set{\mu(T, H^\lam \bot S)}{\lam \in \Gamma_{k,\ell}}$ spans a $\begin{pmatrix}n+\ell\\n \end{pmatrix}$-dimensional $\Q$-space, and the set $\set{\mu(T, H^\mu\bot H^{{\bf0},k-n} \bot S)}{\mu \in \Gamma_{n,\ell}}$ forms a basis.
\end{thm}

\begin{rem}{\rm  (1) In the above theorem, the numbers $a_\tau$ are rational. It follows from the fact that the matrix $\left(q^{\pair{\lam}{\tau}}\right)_{\lam, \tau \in \Lam_{\ell, k}}$ is non-degenerate. (cf. \cite[Proposition 2.7]{BHS}). \\
(2) For symmetric forms, local densities appear in Fourier coefficients of genus theta series. Hence we could pull out the linear indepedence within Siegel Eisenstein series from the linear independence of local densities for symmetric forms (cf. [BHS, \S 4-6]). The author expects that the present results will be helpful for global studies of hermitian modular forms. \\
(3) Local densities and Siegel singular series have a close connection and the latter can be expressed by certain spherical functions (generating functions of local densities). It is possible to obtain the functional equation of Siegel series from such a point of view (cf. symmetric case \cite{HS}, hermitian case \cite{H-Oda}). 
Siegel series are important arithmetic objects and have been well studied (e.g., [Kat]).  The existence of the functional equation is  a rather classical result ([Kar]) and its explicit formula is given for all quadratic extensions over $\Q_p$ in [Ik].   
}
\end{rem}

\slit
Recall the character $\psi$ on $k$ defined in \eqref{character}. For a function $f$ on $V_n(D)$,  we define its Fourier transform if the following integral is well-defined: 
\begin{eqnarray} \label{F-trans def}
\calF(f)(y) = \int_{V_n(D)}\, f(x)\psi(-\pair{x}{y})dx. 
\end{eqnarray}

\mslit
For $A \in V_m$ and $C \in V_n$, we define the Gauss sum by
\begin{eqnarray} \label{G-sum}
\calG(A, C) = \int_{M_{mn}(\calO)}\, \psi(\pair{A[v]}{C})dv, 
\end{eqnarray}
which becomes a rational number as is calculated explicitly in Proposition~\ref{prop: 3.2}. %
Then, by Lemma~\ref{lem 1.2}, we see 
\begin{eqnarray} 
&&
\calF(\calG(A,\;))(B) = q^{-n(n-1)}\mu(B, A), \nonumber\\
&&\label{F-tr Gauss}
\quad \mbox{if } \int_{V_n(D)}\calG(A, x)\,\psi(-\pair{x}{B})dx \mbox{ is well-defined}.
\end{eqnarray}

\mslit
For $x \in V_n$, we define a constant $\nu[x]$ as follows.  
If $x \in V_n(\calO)$, we set $\nu[x]=1$; if $x$ has eigenvalues of negative $\varPi$-exponents and $\varPi^{-\tau_1}, \ldots, \varPi^{-\tau_r}$ are the all, we set 
\begin{eqnarray} \label{def: nu[x]}
\nu[x] = q^{\abs{\tau}}, \quad \abs{\tau} = \sum_{i=1}^r \tau_i  \; (>0).
\end{eqnarray}

\begin{prop} \label{prop: estimate Gauss}
For $A \in X_m$ and $C \in V_n$, it holds
\begin{eqnarray}
\abs{\calG(A, C)} \leq c(A)\, \nu[C]^{-m},
\end{eqnarray}
where, taking $\alp$ as the largest exponent of $\varPi$ among the eigenvalues of $A$,
\begin{eqnarray}
c(A) = q^{mn(\alp+1)}.
\end{eqnarray}
\end{prop}
We admit Proposition~\ref{prop: estimate Gauss} for a while, which will be proved  in \S3.

\begin{prop} \label{prop: F-tr of density}
For $A \in X_m$ and $T \in V_n$, the Fourier transform $\calF(\calG(A, \;))(T)$ is absolutely convergent if $m$ and $n$ satisfy the condition $m \geq 8n-1$. 
Then, for $B \in X_n$, one has $\calF(\calG(A,\; ))(B) = q^{-n(n-1)}\mu(B, A)$.
\end{prop}
\proof
By Proposition~\ref{prop: estimate Gauss}, the right hand side of \eqref{F-tr Gauss} is bounded from the upper, except the constant $c(A)$ by 
\begin{eqnarray} \label{upper bdd}
\int_{V_n}\, \nu[y]^{-m}dy.
\end{eqnarray}
By the map $\vphi_n$, we embed $V_n(D)$ into $M_{2n}(k'), \; k'=k(\eps)$  (\eqref{vphi-n}). 
We recall  
\begin{eqnarray*}
\vphi(\varPi^{-2\ell}) = \twomatrix{\pi^{-\ell}}{0}{0}{\pi^{-\ell}}, \quad
\vphi(\varPi^{-2\ell+1}) = \twomatrix{0}{\pi^{-\ell+1}}{\pi^{-\ell}}{0}, \quad (\ell \geq 1).
\end{eqnarray*}
Hence, if $x \in V_n(D)$ has eigenvalues of negative $\varPi$-exponents, $\vphi_n(x)$ does so, and  the sum of negative $\pi$-exponents of $\vphi_n(x)$ is equal to $\abs{\tau}$, which is the same as in \eqref{def: nu[x]}. 
In general, for $y \in M_{2n}(k')$, one can define $\nu_{k'}(y)$ as follows: if $y$ has eigenvalues of negative $\pi$-exponents and those sum is $c(y)$, then we set $\nu_{k'}(y) = q^{2\abs{c(y)}} $; while we set $\nu_{k'}(y) = 1$ otherwise.
It is known (e.g. \cite[3.14]{Shimura}) that
\begin{eqnarray*}
\int_{M_{2n}(k')} \nu_{k'}(x)^{-s}dx = \prod_{i=1}^{2n}\frac{1-q^{2(-s+i-1)}}{1-q^{2(-s+2n+i-1)}}, \quad ({\rm Re}(s) > 4n-1).
\end{eqnarray*}
Since $\nu[x]^2 = \nu_{k'}(\vphi_n(x))$, we see \eqref{upper bdd} is 
is absolutely convergent if $m \geq 8n-1$, hence in the same range it holds $\calF(\calG(A,\; ))(B) = \mu(B, A)$.
\qed

\begin{prop} \label{prop4}
Assume that $m$ and $n$ satisfy $m \geq 8n-1$, and let $A_i \in X_m(\calO), 1 \leq i \leq N$. Then the following are equivalent:

\noindent
{\rm (i)} As functions of $T$ on $X_n^+$, $\mu(T, A_i), \; 1 \leq i \leq N$ are linearly independent over $\Q$.

\noindent
{\rm (ii)} As functions of $X$ on $V_n$, $\calG(A_i, X), \; 1 \leq i \leq N$ are linearly independent over $\Q$.
\end{prop}

\proof (The strategy is the same as in [BHS], and we
know that $\mu(T, A_i)$ and $\calG(A_i, X)$ are rational numbers.)
It is clear that (i) implies (ii).
Conversely we assume that $\calG(A_i, X)$ are linearly independent over $\Q$ as functions of $X \in V_n$. Under the assumption of $m$ and $n$, $\calG(A_i, X)$ are integralable on $V_n$ and give $\mu(T, A_i)$ at $T \in X_n$, hence $\mu(T, A_i)$ are linearly independent as functions of $T \in V_n$.  In general, $\calG(A, Y) = 1$ for any $A \in X_m^+$ and $Y \in V_n(\pi, 0)$, hence we have
\begin{eqnarray*}
\mu(T, A_i) &=& \int_{V_n}\calG(A_i, X+Y) \psi(-\pair{T}{X})dX 
= \int_{V_n} \calG(A_i, X) \psi(-\pair{T}{X-Y})dX\\
&=& \psi(\pair{T}{Y}) \mu(T, A_i).
\end{eqnarray*} 
On the other hand,  if $T \notin V_n(\calO)$, there exists some $Y \in V_n(\pi, 0)$ for which $\psi(\pair{T}{Y}) \ne 1$; at that time by the above equation we see that $\mu(T, A_i) = 0$. Since the set $V_n(\calO) \cap X_n$ is dense in $V_n(\calO)$, and Fourier transform of integralable function is continuous, we see that $\set{\mu(T, A_i)}{1 \leq i \leq N}$ are linearly independent as functions on $X_n^+ = X_n \cap V_n(\calO)$. 
\qed

\vspace{2cm}
\Section{Estimate of Gauss sums and the proof of Theorem~\ref{th: lin indep}}
Let $A \in V_m$ and $C \in V_n$. We see that $\calG(A, C) = \calG(C, A)$ is determined by the $GL_m(\calO)$-orbit containing $A$ and the $GL_n(\calO)$-orbit containing $C$. Further it is clear that $\calG(A, C)$ decomposes into products if $A$ or $C$ decomposes into orthogonal sums, i.e., 
\begin{eqnarray} \label{decom G-sum}
\calG(\bot_{i=1}^r A_i, \bot_{j=1}^t C_j) = \prod_{i=1}^r\prod_{j=1}^t\, \calG(A_i, C_j).
\end{eqnarray} 
Thus we may assume that $A \in V_m$ and $C \in V_n$ have the following shape:
\begin{eqnarray}
&&
A = \gen{\varPi^{\alp_1}}\bot \cdots \bot \gen{\varPi^{\alp_r}} \bot 
\twomatrix{}{\varPi^{\alp_{r+2}}}{-\varPi^{\alp_{r+2}}}{} \bot \cdots \bot \twomatrix{}{\varPi^{\alp_{r+2s}}}{-\varPi^{\alp_{r+2s}}}{},\nonumber \\
&&
\quad (r\geq 0, s \geq 0, r+2s = m, \; 2\vert \alp_\ell \mbox{ for } 1 \leq \ell \leq r, \; 2\not{\mid} \alp_{r+2\mu}=\alp_{r+2\mu-1} \mbox{ for }1 \leq \mu \leq s) \nonumber \\
&&
C = \gen{\varPi^{\beta_1}}\bot\cdots\bot\gen{\varPi^{\beta_t}} \bot \twomatrix{}{\varPi^{\beta_{t+2}}}{-\varPi^{\beta_{t+2}}}{} \bot \cdots \bot\twomatrix{}{\varPi^{\beta_{t+2v}}}{-\varPi^{\beta_{t+2v}}}{},\nonumber \\
&& \label{D-form1} 
\quad (t\geq 0, v \geq 0, t+2v = m, \; 2\vert \beta_i \mbox{ for } 1 \leq i \leq t, \; 2\not{\mid} \beta_{t+2j}=\beta_{t+2j-1} \mbox{ for }1 \leq j \leq v).\nonumber\\
\end{eqnarray}
If some elementary divisor of $A$ or $C$ is $0$, we understand $0 = \varPi^\infty$ in the diagonal component in \eqref{D-form1}, i.e., some of $\varPi^{\alp_i}=0$ or $\varPi^{\beta_i}=0$. 
We set
\begin{eqnarray} \label{I for D}
I(a) &=& \int_{\calO} \psi(\varPi^a \Nrd(x))dx,\quad (a \in 2\Z) \nonumber\\
\label{J for D}
J(b) &=& \int_{\calO\times\calO}\, \psi(\Trd(\varPi^b xy))dxdy
= \int_{\calO\times\calO}\, \psi(\Trd(xy\varPi^b))dxdy, \quad (b \in \Z),
\end{eqnarray}
then we have the following.

\begin{lem} \label{lem: 3.1}
Assume that $a \in \Z$ is even and $b \in \Z$. Then
\begin{eqnarray*}
I(a) = (-q)^{\min\{0, \, a+1\}}, \quad J(b) = q^{2\min\{0, \, b+1\}}.
\end{eqnarray*}
\end{lem}

\proof
First we note that $\Nrd(\calO^\times) = \calo^\times$ and $\Trd(\calP^m) = \frp^{[(m+1)/2]}$ for $m \in \Z$. 
Then, it is clear that $I(a) = 1$ if $a \geq 0$ and $J(b) = 1$ if $b \geq -1$.

\mmslit
Assume that $a = -2\ell \leq -2$. Then
\begin{eqnarray*}
I(a) &=& \int_{\calO} \psi(\pi^{-\ell}\Nrd(x))dx 
= q^{-4\ell}\sum_{\ol{x} \in \calO/\calP^{2\ell}}\, \psi(\pi^{-\ell}\Nrd(x)) \nonumber \\
&=& \label{I-D1}
q^{-4\ell}\sum_{k=0}^{\ell-1}\sum_{\ol{u} \in (\calO/\calP^{2\ell-k})^\times}\, \psi(\pi^{-\ell} \Nrd(\varPi^k u) + vol(\calP^\ell/\calP^{2\ell}) \nonumber\\
&=&
q^{-4\ell}\sum_{k=0}^{\ell-1}q^{3\ell-k} (1+q^{-1})\underbrace{\sum_{\eps \in (\calo/\frp^{\ell-k})^\times}\, \psi(\chi^{-\ell+k}\eps)}_{(*1)} + q^{-2\ell}.
\end{eqnarray*}
Since $(*1) = 0$ unless $k=\ell-1$ and $(*1) = -1$ if $k=\ell-1$, we obtain
\begin{eqnarray*}
I(a) &=&
q^{-4\ell}(-q^{2\ell+1}(1+q^{-1}))+q^{-2\ell} = -q^{-2\ell+1} = -q^{a+1}.
\end{eqnarray*}
Next we consider $J(b)$ with $b = -\ell \leq -2$.
\begin{eqnarray*}
J(b) 
&=&
q^{-4(\ell-1)} \sum_{\ol{x}, \ol{y} \in \calO/\calP^{\ell-1}}\, \psi(\Trd(\varPi^{-\ell}xy) \nonumber\\
&=& \label{J-1}
q^{-4(\ell-1)} \sum_{\ol{x} \in \calO/\calP^{\ell-1}}\, \left(
\sum_{k=0}^{\ell-2}\sum_{\ol{u} \in (\calO/\calP^{\ell-k-1})^\times}\, \psi(\Trd(\varPi^{-\ell+k}xu)) + 1 \right)\\
&=& \label{J-2}
q^{-4(\ell-1)} \sum_{k=0}^{\ell-2}\, q^{2(\ell-k-1)}(1-q^{-2}) \underbrace{\sum_{\ol{x} \in \calO/\calP^{\ell-1}}\, \psi(\Trd(\varPi^{-\ell+k}x))}_{(*2)} + q^{-2(\ell-1)}.
\end{eqnarray*}
Since we have the following additive epimorphism
\begin{eqnarray*}
\calO/\calP^{\ell-1} \longrightarrow \frp^{-m}/\calo, \; x + \calP^{\ell-1} \longmapsto \Trd(\varPi^{-\ell+k}x), \quad m = \left[\frac{\ell-k-1}{2}\right] >0,
\end{eqnarray*}
and $\psi$ is nontrivial on $\frp^{-m}/\calo$, we see $(*2) = 0$ for $0 \leq k \leq \ell-2$.
Thus, we obtain $J(b) = q^{-2(\ell-1)} = q^{2(b+1)}$.
\qed

\begin{prop} \label{prop: 3.2}
Assume that $A$ and $C$ are given as in the shape of \eqref{D-form1}. Then
\begin{eqnarray*} \label{G-sumD1}
\calG(A, C) &=&
(-1)^c \prod_{i=1}^m \, \prod_{j=1}^n\, q^{\min\{0, \alp_i+\beta_j+1\}},
\end{eqnarray*}
where 
\begin{eqnarray*} \label{const c}
c = \sharp\set{(\ell, i)}{\alp_\ell + \beta_i <-1, \; 1 \leq \ell \leq r, \; 1\leq i \leq t}.
\end{eqnarray*}
\end{prop}

\proof
For $x = a+b\eps + c\varPi+d\varPi\eps \in D$ with $a, b, c, d \in k$, we set $x' = a-b\eps + c\varPi -d\varPi\eps$. Then, $\varPi x = x' \varPi$, ${x'}^* = {x^*}'$, and $x'$ runs over $\calO$ when $x$ does so.

\mmslit
Assume that $a$ and $b$ are even and $r$ and $s$ are odd, and recall $h^r$ of \eqref{H-lam}.  
Clearly $\calG(\varPi^a, \varPi^b) = I(a+b)$. Further, we have
\begin{eqnarray*}
\calG(h^r, \varPi^a) &=& \calG(\varPi^a, h^r) = \int_{\calO\times\calO}^, \psi(\Trd(u^*\varPi^rv\varPi^a)) dudv\\
&=&
J(a+r),\\
\calG(h^r, h^s) &=& \int_{m_2(\calO)}\, \psi(\pair{h^r[u]}{h^s})du \quad (u = \twomatrix{u_1}{u_2}{v_1}{v_2})\\
&=&
\int_{M_2(\calO)}\, \psi(\Trd(-(u_1^*\varPi^rv_2-v_1^*\varPi^ru_2)\varPi^s))du\\
&=&
J(r+s)^2.
\end{eqnarray*}
By these calculation and Lemma~\ref{lem: 3.1}, we obtain
\begin{eqnarray*}
\calG(A, C) &=&
\prod_{i=1}^t  \left(\prod_{\ell=1}^r I(\alp_\ell+\beta_i) \prod_{\mu=1}^s J(\alp_{r+2\mu}+\beta_i) \right)\nonumber\\
&& \times
\prod_{j=1}^v\left(\prod_{\ell=1}^r J(\alp_\ell+\beta_{t+2j}) \prod_{\mu=1}^s J(\alp_{r+2\mu}+\beta_{t+2j})^2 \right)\\
&=&
(-1)^c\prod_{i=1}^t \left(\prod_{\ell=1}^r q^{\min\{0, \alp_\ell+\beta_i+1\}}\prod_{\mu=1}^s q^{\min\{0, 2(\alp_{r+2\mu}+\beta_i+1)\}}\right) \nonumber \\
&& 
\times
\prod_{j=1}^v\left(\prod_{\ell=1}^r q^{\min\{0, 2(\alp_\ell+\beta_{t+2j}+1)\}} \prod_{\mu=1}^s q^{2\min\{0, 2(\alp_{r+2\mu}+\beta_{t+2j}+1)\} } \right) \nonumber\\
&=& \label{G-sum2}
(-1)^c \prod_{i=1}^m \, \prod_{j=1}^n\, q^{\min\{0, \alp_i+\beta_j+1\}}.
\end{eqnarray*}
\qed

\slit
{\it Proof of Proposition~2.2}.
Since $A \in X_m$, we may take $\alp = \max\set{\alp_i}{1 \leq i \leq m}$, and 
\begin{eqnarray*}
&& 
q^{\min\{0, \alp_i+\beta_j+1\}} = \left\{\begin{array}{ll}
q^{\alp_i+\beta_j+1} & \mbox{if } \alp_i+\beta_j +1 < 0\\
1 & \mbox{if } \alp_i+\beta_j +1 \geq 0
\end{array} \right\} \leq 1, \nonumber \\[2mm]
&&\label{estimate1}
q^{\min\{0, \alp_i+\beta_j+1\}} 
\leq
q^{\alp_i+\beta_j+1}  \leq  q^{\alp+\beta_j+1}.
\end{eqnarray*}
Hence by Proposition~\ref{prop: 3.2}, we have
\begin{eqnarray*}
\abs{\calG(A, C)} &\leq& \prod_{i=1}^m\, 
\prod_{{\scriptsize \begin{array}{l}1\leq j \leq n\\ \beta_j < 0 \end{array}}}\,
q^{\alp+\beta_j+1}  \nonumber\\
&\leq& \label{absGauss-D}
\prod_{{\scriptsize \begin{array}{l}1\leq j \leq n\\ \beta_j < 0 \end{array}}}\,
q^{m(\alp+\beta_j+1)}  \leq c(A)\, \nu[C]^{-m},
\end{eqnarray*}
where $c(A)= q^{mn(\alp+1)}$. 
\qed

\slit
We prepare some more notation.
For each $X \in V_n \backslash V_n(\calO)$, let $\varPi^{-\tau_1}, \ldots, \varPi^{-\tau_r}$ be those elementary divisors of $X$ which are of negative powers of $\varPi$ and $\tau_1\geq \cdots \tau_r \geq 1$. Set 
\begin{eqnarray} \label{tau(X)}
\tau = \tau(X) = \left\{\begin{array}{ll}
0 \in \Lam_n^+ & \mbox{if} X \in V_n(\calO)\\
(\tau_1, \ldots, \tau_r) \in \Lam_r^+(\subset \Lam_n^+) & \mbox{if} X \notin V_n(\calO), 
\end{array}\right.
\end{eqnarray}
and 
\begin{eqnarray} \label{sigma(X)}
\sigma = \sigma(X) = \left\{\begin{array}{ll}
0 \in \Lam_n^+ & \mbox{if} X \in V_n(\calO)\\
(\tau_1-1, \ldots, \tau_r-1) \in \Lam_r^+(\subset \Lam_n^+) & \mbox{if} X \notin V_n(\calO). 
\end{array}\right.
\end{eqnarray}
We note that $\nu[X] = q^{\abs{\tau}}$ (cf. \eqref{def: nu[x]}). 

\mslit
We recall $H^\lam$ for $\lam \in \Gamma_k$ (cf. \eqref{H-lam}) and assume that the elementary divisors of $X \in V_n$ are given by $\{ \beta_j \}$ of the shape $C$ in \eqref{D-form1}. Then, by Proposition~\ref{prop: 3.2}, we have
\begin{eqnarray*}
\calG(H^\lam, X)&=&
\prod_{i=1}^k \prod_{j=1}^n q^{2\min\{0, \lam_i +\beta_j+1\}}
=
q^{\sum_{t \geq 0}\, m_\lam(t) \sum_{j=1}^n\, 2\min\{0, t+\beta_j+1\}}
 = q^{c(\lam, X)},
\end{eqnarray*}
where 
\begin{eqnarray*}
c(\lam, X) &=& \sum_{t \geq 0}\, m_\lam(t) \sum_{j=1}^n\, 2\min\{0, t+\beta_j+1\}.
\end{eqnarray*}
If $X \in V(\calO)$, clearly $c(\lam, X) = 0$. When $X \notin V(\calO)$, we have
\begin{eqnarray}
c(\lam, X) &=&
2\sum_{t \geq 0}\, m_\lam(t) \sum_{\scriptsize{\begin{array}{l} 1 \leq j \leq n \\ t + 1 + \beta_j < 0 \end{array}}} (t+1-\tau_j) \qquad (m_\lam(t) = \sharp\set{i}{\lam_i = t}) \nonumber \\
&=&
2\sum_{t \geq 0}\, m_\lam(t) \sum_{\scriptsize{\begin{array}{l} 1 \leq i \leq r,\\ \sigma_i > t \end{array}}}\, (t-\sigma_i)  \label{temp1} .
\end{eqnarray}
By a combinatorial calculation for \eqref{temp1}, which is the same as in \cite[p.58]{BHS}, we obtain
\begin{eqnarray} \label{G(lam, X)}
\calG(H^\lam, X) = q^{2\pair{\what{\lam}}{\what{\sigma}}-2k\abs{\sigma}},
\end{eqnarray}
where the last expression is valid even if $X \in V(\calO)$, since $\sigma = 0$ in that case.

\mslit
{\it Proof of Theorem~\ref{th: lin indep}}\\
 (1) Under the given condition on $k, n$ and $r$, the condition of Proposition~\ref{prop: F-tr of density} is satisfied as $m=2k+r \geq 8n-1$, hence local densities $\mu(T, H^\lam \bot S)$ are expressed as the Fourier transform of Gauss sums (cf. \eqref{F-tr Gauss}).
By \eqref{G(lam, X)}, we have
\begin{eqnarray*} 
q^{-n(n-1)} \mu(T, H^\lam \bot S) 
&=& \int_{V_n}q^{2\pair{\what{\lam}}{\what{\sigma(X)}}}q^{-2k\abs{\sigma(X)}} \calG(S, X) \psi(-\pair{T}{X})dX ,
\end{eqnarray*} 
where  $\pair{\what{\lam}}{\what{\sigma(X)}}$ is determined by $\what{\sigma(X)}_i, \;
 1 \leq i \leq \ell$ with $\ell$ coming from $\lam \in \Gamma_{k,\ell}$. 
Hence we have 
\begin{eqnarray*} 
\lefteqn{\mu(T, H^\lam \bot S)}\\[2mm]
  && =
q^{n(n-1)}
\sum_{\tau \in \Gamma_{n, \ell}} q^{2\pair{\what{\lam}}{\what{\tau}}}\underbrace{\int_{\set{X\in V_n}{\what{\sigma(X)}_i = \tau_i, \; 1 \leq i \leq \ell}}\, q^{-2k\abs{\sigma(X)}} \calG(S, X) \psi(-\pair{T}{X})dX}_{(*)}.
\end{eqnarray*}
The value $(*)$ is independent of $\lam$ and only depends on $k, \tau, S$ and $T$. Putting $q^{n(n-1)} \cdot (*)$ as $a_\tau$ we have the required formulation.
 
\mslit
(2)   Denote by $W$ the $\Q$-space spanned by functions $\set{\mu(T, H^\lam\bot S)}{\lam \in \Gamma_{k,\ell}}$ on $T \in X_n$. Then, by Proposition~\ref{prop4}, $W$ is isomorphic to the $\Q$-space spanned by functions $\set{\calG(H^\lam\bot S, X)}{\lam \in \Gamma_{k,\ell}}$ on $X \in V_n$. 
On the other hand, the space $W$ is isomorphic to the $\Q$-space $W_0$ spanned by functions $\set{ q^{2\pair{\what{\lam}}{\what{\tau}}}}{\lam \in \Gamma_{k,\ell}}$ on $\tau \in \Gamma_{n,\ell}$ by (1).
The space $W_0$ coincides with that considered in [BHS, p.59], where $q^2$ should be replaced by $q$, and it is proved there that $\dim(W_0) = \sharp(\Gamma_{n,\ell}) = \begin{pmatrix}n+\ell\\ n \end{pmatrix}$. Returning to the space $W$, we see that the set $\set{\mu(T, H^\mu\bot H^{{\bf0}, k-n}\bot S)}{\mu \in \Gamma_{n,\ell}}$ 
forms a basis.
\qed
\vspace{2cm}
\Section{Kitaoka series}
Let us consider the following formal power series for $A \in X_m$ and $B \in X_n$ with $m \geq n$.
\begin{eqnarray}
P(B, A; X) = \sum_{r \geq 0}\, \mu(\pi^rB, A)X^r.
\end{eqnarray}
Kitaoka introduced a similar power series for symmetric forms (\cite{Ki}), conjectured it is rational and proved for a special case, hence we call $P(B, A; X)$ as Kitaoka series. Then B\"ocherer and Sato proved the rationality by using Denif's theory and calculated the denominators for certain cases (\cite{BS}). And the author determined the denominators in \cite{Kitaoka} and \cite{Kitaoka2} by an elementary method, where the rationality is also assured. Here, following the same method, we prove the rationality and determine the denominators for quaternion hermitian case.  
Our result is the following.

\begin{thm} \label{th: KitaokaS}
Assume that $A \in X_m$ and $B \in X_n$ with $m \geq n$. Then 
\begin{eqnarray*}
P(B, A; X) \times \prod_{i=0}^n\, (1 - q^{(n-i)(2n+2i-2m-1)}X)
\end{eqnarray*}
is a polynomial in $X$.
\end{thm}

\mmslit
Because of \eqref{shift}, it suffices to consider integral forms.
We recall the additive character $\psi_\ell$ on $k$ of conductor $\frp^\ell$ which induces the character $\chi_\ell$ on $\calo/\frp^\ell$ (cf. \eqref{character}). 
For a locally constant compactly supported function $f$ on $V_n$ (i.e., for $f \in \calS(V_n)$), we define its Fourier transform with respect to $\psi_\ell$ by
\begin{eqnarray}
(f)^\wedge_\ell(z) = \int_{V_n}\, f(y)\psi_\ell(-\pair{y}{z})dy,
\end{eqnarray}
where $dy$ is the Haar measure on $V_n$ normalized by $vol(V_n(\calO)) = 1$. 

\mslit
For $A \in X_m^+$ and $C \in V_n(\calO)$, we define (finite Gauss sum) 
\begin{eqnarray} \label{finite G-sum}
\calS_\ell(A, C) = \sum_{\ol{v} \in M_{mn}(\calO/\calP^{2\ell})}\, \chi_\ell(\pair{A[v]}{C}).
\end{eqnarray}%
This is well-defined, since $\chi_\ell(\pair{A[v]}{C})=\psi_\ell(\pair{A[v]}{C})$ is stable on the cosets of \\$M_{mn}(\calO)/M_{mn}(\calP^{2\ell}) \cong M_{mn}(\calO/\calP^{2\ell})$. Further, since $\psi_\ell(\pair{A[v]}{C}) = \psi(\pi^{-\ell}\pair{A[v]}{C}) = \psi(\pair{\pi^{-\ell}A[v]}{C})$,  we see 
\begin{eqnarray}
\calS_\ell(A, C) = q^{4\ell mn} \calG(\pi^{-\ell}A, C), \; (A \in X_m^+, C \in V_n(\calO)),
\end{eqnarray}
where $\calG(\;,\;)$ is defined in \eqref{G-sum}. Thus, by Proposition~\ref{prop: 3.2}, we have the following.

\begin{prop} \label{prop: calS(A,C)}
For $A = \pi^\alp \in X_m^+, \, (\alp\in \Lam_m^+)$ and $C = \pi^\beta \in X_n^+, \, (\beta \in \Lam_n^+)$  {\rm (}cf. \eqref{pi-gamma}{\rm)}, it holds
\begin{eqnarray*} \label{S(A,C)-value}
&&
\calS_\ell(A, C) = q^{2\ell mn} \cdot (-1)^c\,  \prod_{i=1}^m\, \prod_{j=1}^n\, q^{\min\{2\ell, \, \alp_i+\beta_j+1\}}, \\
&&
\qquad c = \sharp\set{(i,j)}{\mbox{$\alp_i$ and $\beta_j$ are even, and } 2\ell > \alp_i+\beta_j +1}. \nonumber 
\end{eqnarray*}
If  $2\ell \geq \alp_1+\beta_1+1$, then
\begin{eqnarray*} \label{S(A,C)-v2}
\calS_\ell(A, C) = q^{2\ell mn} \cdot (-1)^c\, q^{n\abs{\alp}+m\abs{\beta}+mn}.
\end{eqnarray*}
\end{prop}

\slit
The group $GL_n(\calO/\calP^{2\ell})$ acts on $V_n(\calO/\calP^{2\ell})$, and the coset space bijectively corresponds to the set 
\begin{eqnarray} \label{Lam-n,2l}
\Lam_{n, 2\ell}^+ = \set{\gamma \in \Lam_n^+}{\gamma_1 \leq 2\ell}.
\end{eqnarray}
For $\gamma \in \Lam_{n,2\ell}^+$, we define 
\begin{eqnarray}
N_\ell^{pr}(\pi^\gamma, \pi^\gamma) = \sharp\set{g \in GL_n(\calO/\calP^{2\ell})}{\pi^\gamma[g]-\pi^\gamma \in V_n(\pi, \ell)}.
\end{eqnarray}
We note here that $N_\ell^{pr}(\pi^\gamma, \pi^\gamma) = N_\ell(\pi^\gamma, \pi^\gamma)$ if $\gamma_1 < 2\ell$ (cf. \eqref{loc-dens3}).
We will rewrite the local density as follows.

\begin{prop} \label{prop: K-3}
Assume that $A \in X_m^+$ and $B \in X_n^+$ with $m \geq n$ and that $\ell$ satisfies the condition \eqref{assump}. Then one has
\begin{eqnarray*} 
\mu(B, A) &=&    \label{muBA}
 q^{-\ell n(4m-2n+1)-n(n-1)} N_\ell(B, B) \sum_{\gamma \in {\Lam^+_{n,2\ell}}}\, \frac{\calS_\ell(A, \pi^\gamma)}{N^{pr}_\ell(\pi^\gamma,\pi^\gamma)}\, (ch_B)^\wedge_\ell(\pi^\gamma),
\end{eqnarray*}
where $ch_B$ is the characteristic function of $K_n\cdot B$.
\end{prop}

\proof
Since $\ell$ is big enough, we have by \eqref{N in finite}
\begin{eqnarray*}
\mu(B, A) &=& q^{\ell n(4m-2n+1)-n(n-1)}N_\ell(B, A)\nonumber\\
&=&  
q^{-4\ell mn -n(n-1)} \sum_{\ol{y} \in V_n(\calO/\calP^{2\ell})}\, \calS_\ell(A, y) \, \chi_\ell(\pair{-B}{y}) \nonumber\\
&=&
q^{-4\ell mn -n(n-1)} \sum_{\gamma \in \Lam_{n, 2\ell}^+}\, \calS_\ell(A, \pi^\gamma)
\underbrace{ \sum_{\scriptsize{\begin{array}{c} \ol{y}\in V_n(\calO/\calP^{2\ell})\\ y \sim \pi^\gamma\end{array}}}\, \chi_\ell(-\pair{y}{B})}_{(*1)},
\end{eqnarray*}
where  $y \sim \pi^\gamma$ means that $\ol{y} \in GL_n(\calO/\calP^{2\ell}) \cdot \pi^\gamma, \; (\gamma \in \Lam_{n, 2\ell}^+$).
Here 
\begin{eqnarray*}
(*1) &=& 
\frac{1}{N_\ell^{pr}(\pi^\gamma, \pi^\gamma)} \sum_{\ol{g} \in GL_n(\calO/\calP^{2\ell})}\, \chi_\ell(-\pair{\pi^\gamma}{g^*Bg}) \nonumber\\
&=&
\frac{N_\ell^{pr}(B, B)}{N_\ell^{pr}(\pi^\gamma, \pi^\gamma)}\, \sum_{\scriptsize{\begin{array}{c} \ol{z}\in V_n(\calO/\calP^{2\ell})\\z \sim B \end{array}}}\, \chi_\ell(-\pair{\pi^\gamma}{z}). 
\end{eqnarray*}
Thus we obtain
\begin{eqnarray} 
\lefteqn{\mu(B, A) }\nonumber\\
&=& \label{mu(B,A)-2}
q^{-4\ell mn-n(n-1)} N_\ell^{pr}(B, B) \sum_{\gamma \in \Lam_{n,2\ell}^+} \frac{S_\ell(A,\pi^\gamma)}{N_\ell^{pr}(\pi^\gamma, \pi^\gamma)}\, \sum_{\scriptsize{\begin{array}{c} \ol{z}\in V_n(\calO/\calP^{2\ell})\\z \sim B\end{array}}}\, \chi_\ell(\pair{\pi^\gamma}{z}).
\end{eqnarray}
By the assumption \eqref{assump} for $\ell$, the $K_n$-orbit
$K_n\cdot B$ is decomposed as a finite union of $V_n(\pi, \ell)$-cosets, and 
\begin{eqnarray}\label{F-trans}
(ch_B)^\wedge_\ell (\pi^\gamma) = \int_{K_n \cdot B}\psi_\ell(-\pair{z}{\pi^\gamma})dz = q^{-\ell n(2n-1) } \sum_{\scriptsize{\begin{array}{c} \ol{z}\in V_n(\calO/\calP^{2\ell})\\z \sim B \end{array}}}\, \chi_\ell(\pair{\pi^\gamma}{z}).
\end{eqnarray}
By \eqref{mu(B,A)-2} and \eqref{F-trans}, we obtain the required formula for $\mu(B, A)$.
\qed

\slit
We note the following results of local densities (cf. \cite[Th.2.3, Prop.2.4, Prop.2.6]{Q}). 

\begin{lem}\label{lem: N-l} 
{\rm (1)} Assume that $A \in V_m(\calO)$, $B \in V_n(\calO)$ with $m \geq n$ and $e \in \N$. Then, 
\begin{eqnarray*}
N_{\ell+e}(\pi^eB, \pi^eA) = q^{4emn} N_\ell(B, A).
\end{eqnarray*}
{\rm (2)} Assume that $\alp \in \Lam_n^+$ has the shape $\alp = \beta\gamma$ with $\beta \in \Lam_i^+$, $\gamma \in \Lam_{n-i}^+$ and $\beta_i > \gamma_1$. Then,
\begin{align}
&
N_\ell^{pr}(\pi^\alp, \pi^\alp) = q^{(4\ell+2) i(n-i)+2i\abs{\gamma}}\, N_\ell^{pr}(\pi^\beta, \pi^\beta) N_\ell^{pr}(\pi^\gamma, \pi^\gamma). \nonumber
\end{align}
{\rm (3)} Assume that $\alp \in \Lam_{n, \ell}^+$ has the shape
\begin{eqnarray*}
\alp = (\alp_1^{m_1} \cdots \alp_s^{m_s}),\;  \alp_1>\cdots >\alp_s, \; m_i \geq 1 \; \any i.
\end{eqnarray*}
Then
\begin{eqnarray*}
N_{\ell}^{pr}(\pi^\alp, \pi^\alp) &=&
q^{m_\alp} 
\times \prod_{i=2}^r\, \left\{\begin{array}{ll}
w_{m_i}(-q^{-1}) & \mbox{if } 2 \mid \alp_i\\
w_{m_i/2}(q^{-4}) & \mbox{if } 2 \not{\mid} \alp_i
\end{array}\right\} \nonumber\\
&&
\times
\left\{\begin{array}{ll}
w_{m_1}(-q^{-1}) & \mbox{if } 2 \mid \alp_1 < 2\ell\\
w_{m_1/2}(q^{-4}) & \mbox{if } 2 \not{\mid} \alp_1 < 2\ell\\
 w_{m_1}(q^{-2}) & \mbox{if } \alp_1 = 2\ell
\end{array}\right\},
\end{eqnarray*}
where $w_m(t) = \prod_{i=1}^m\, (1 - t^i)$  and
\begin{eqnarray*}
m_\alp = \ell n(2n+1)+n(n-1) + 2\sum_{i=1}^n\, (i-1)\alp_i +\frac12\abs{\alp}+\frac12\sharp\set{i}{\alp_i \mbox{ is odd} }.
\end{eqnarray*}
\end{lem}
For simplicity, we introduce the following symbol 
\begin{eqnarray} \label{H-k}
H_k = \twomatrix{0}{\varPi}{-\varPi}{0} \bot \cdots \bot \twomatrix{0}{\varPi}{-\varPi}{0} \in X_{2k}.
\end{eqnarray} 
By Proposition~\ref{prop: calS(A,C)} and Lemma~\ref{lem: N-l}(3), we have the following data, which we use for the calculation of $P(B, A; X)$ in \S 5.

\begin{lem} \label{lem: low ratio S-N}
Assume that $2\ell \geq \alp+2$, where  $\alp$ is the maximal $\varPi$-exponent of all the elementary divisors of $A \in X_m^+$. Then
\begin{eqnarray*}
\dfrac{S_\ell(A, 1_n)}{N_\ell(1_n, 1_n)} = c_0\, q^{\ell n(2m-2n-1)},\qquad
\dfrac{S_\ell(A, H_{n/2})}{N_\ell(H_{n/2}, H_{n/2})} = c_1\, q^{\ell n(2m-2n-1)},
\end{eqnarray*}
where $c_0$ and $c_1$ are constants determined by $A$ and $n$, and independent of $\ell$.  
\end{lem}

\slit
\begin{prop} \label{prop: P(B,A)}
Assume that $\ell$ satisfies the condition \eqref{assump} for $B$. Then
\begin{eqnarray*}
P(B,A; X) &=& 
q^{-\ell n(4m-2n+1)-n(n-1)} N_\ell(B, B) \sum_{\wt{\gamma} \in {\Lam^+_{n,2\ell}}}\, (ch_B)^\wedge_\ell(\pi^{\wt{\gamma}})\, P_{\wt{\gamma}}(A; X),
\end{eqnarray*}
where
\begin{eqnarray}  \label{P-wt-gamma}
P_{\wt{\gamma}}(A; X) &=& \sum_{r=0}^\infty (q^{4n(n-m)}X)^r\, 
\sum_{\gamma \in \Lam(\wt{\gamma},r)}\,\frac{\calS_{\ell+r}(A, \pi^\gamma)}{N_{\ell+r}^{pr}(\pi^\gamma,\pi^\gamma)},\\
\Lam(\wt{\gamma},r) &=& \set{\gamma \in \Lam_{n, 2(\ell+r)}^+}{\gamma \equiv \wt{\gamma} \pmod{2\ell}}. \nonumber
\end{eqnarray}
\end{prop}
\proof
For $\pi^rB$, $\ell + r$ satisfies the condition \eqref{assump}, and we calculate $\mu(\pi^rB, A)$ by using Proposition~\ref{prop: K-3}:
\begin{eqnarray*}
\mu(\pi^rB, A) &=& q^{-(\ell+r)n(4m-2n+1)-n(n-1)} N_{\ell+r}(\pi^rB, \pi^rB)\\
&&\qquad
 \times \sum_{\gamma \in {\Lam^+_{n,2(\ell+r)}}}\, \frac{\calS_{\ell+r}(A, \pi^\gamma)}{N^{pr}_{\ell+r}(\pi^\gamma,\pi^\gamma)}\, (ch_{\pi^rB})^\wedge_{\ell+r}(\pi^\gamma).
\end{eqnarray*}
Here we see
\begin{eqnarray*}
N_{\ell+r}(\pi^rB, \pi^rB) &=& N_{\ell+r}^{pr}(\pi^rB, \pi^rB) = q^{4rn^2}N_\ell(B,B),\; (\mbox{by lemma~\ref{lem: N-l}(1)})\\
(ch_{\pi^rB})^\wedge_{\ell+r}(\pi^\gamma) &=& \int_{V_n}ch_{\pi^rB}(y)\psi_{\ell+r}(-\pair{y}{\pi^\gamma})dy
= q^{-rn(2n-1)}\int_{V_n}ch_B(y)\psi_\ell(-\pair{y}{\pi^\gamma})dy\\
&=&
q^{-rn(2n-1)} (ch_B)^\wedge_\ell(\pi^\gamma),
\end{eqnarray*}
and the last term depends only on $\gamma \pmod{2\ell}$. Thus we have
\begin{eqnarray*}
\mu(\pi^r B, A) &=& q^{-\ell n(4m-2n+1)-n(n-1)} q^{4rn(n-m)} N_\ell(B, B) \nonumber \\
&&   \label{mn-rBA2}
\times  \sum_{\wt{\gamma} \in \Lam_{n, 2\ell}^+}\, (ch_B)^\wedge_\ell(\pi^{\wt{\gamma}})
\sum_{{\scriptsize\begin{array}{l}\gamma \in \Lam_{n, 2(\ell+r)}^+\\ \gamma \equiv \wt{\gamma} \pmod{2\ell} \end{array}}}\,\frac{\calS_{\ell+r}(A, \pi^\gamma)}{N_{\ell+r}^{pr}(\pi^\gamma,\pi^\gamma)}.
\end{eqnarray*}
Hence we obtain the required expression for $P(B, A; X)$ by the definition of $P_{\wt{\gamma}}(A;X)$.
\qed

\slit
Fix a $\wt{\gamma} \in \Lam_{n,2\ell}^+$ and set $i = m_{\wt{\gamma}}(2\ell) \left(=\sharp \set{j}{1 \leq j \leq n, \; \wt{\gamma}_j = 2\ell}\right)$. 
Then we may write $\wt{\gamma} = (2\ell)^iv$ with some $v\in \Lam_{n-i, 2\ell-1}^+$, and
\begin{eqnarray}
\Lam(\wt{\gamma}, r)
&=& \label{gamma-rho}
\left\{\begin{array}{ll}
\{\wt{\gamma}\} = \{v\} & \mbox{if } i=0 \\
\set{((2\ell)^i+\rho)v \in \Lam_{n, 2(\ell+r)}^+}{\rho \in \Lam_{i, 2r}^+} & \mbox{if } i \geq 1,
\end{array}\right.
\end{eqnarray}
We calculate the main quotient term of $P_{\wt{\gamma}}(A; X)$ as follows.

\begin{lem} \label{lem: ratio S-N}
Assume that $\wt{\gamma} \in \Lam_{n,2\ell}^+$ and $i = m_{\wt{\gamma}}(2\ell)$. Take any $\gamma \in \Lam(\wt{\gamma},r)$ written as in \eqref{gamma-rho}. Then, for sufficiently large $r$ with respect to $A$, one has
\begin{eqnarray*}
\dfrac{\calS_{\ell+r}(A, \pi^\gamma)}{N_{\ell+r}^{pr}(\pi^\gamma, \pi^\gamma)}
 = 
\left\{ \begin{array}{ll}
a_0 \, q^{rn(2m-2n-1)} & \mbox{if } i = 0\\
a_i \, q^{r(n-i)(2m-2n-2i-1)} \dfrac{\calS_r(A, \pi^\rho)}{N_r^{pr}(\pi^\rho, \pi^\rho)}  & \mbox{if } 1 \leq i \leq n,
\end{array}\right.
\end{eqnarray*}
where $a_i$ is a constant determined by $\ell, A$ and $\wt{\gamma}$, and independent of $r$. 
\end{lem}

\proof
Assume that $\gamma \in \Lam(\wt{\gamma},r)$ is written as in \eqref{gamma-rho}. Then we have 
\begin{eqnarray*}
\calS_{\ell+r}(A, \pi^\gamma)
& = &
\left\{\begin{array}{ll}
\calS_{\ell+r}(A, \pi^v) & \mbox{if } i=0\\
q^{4\ell mi}\, \calS_r(A, \pi^\rho)\, \calS_{\ell+r}(A, \pi^v) &\mbox{if } 1 \leq i \leq n-1\\
q^{4\ell mn}\, \calS_r(A, \pi^\rho) & \mbox{if } i = n.
\end{array} \right.
\end{eqnarray*}
When $r$ is big enough for the eigenvalues of $A$ and $\ell$, we see by Proposition~\ref{prop: calS(A,C)}
\begin{eqnarray} \label{calS big-r}
\calS_{\ell+r}(A, \pi^\gamma) = \left\{\begin{array}{ll}
c'_0 q^{2rmn} & \mbox{if } i=0\\
c'_i q^{2rm(n-i)}\,\calS_r(A, \pi^\rho) & \mbox{if } 1 \leq i \leq n,
\end{array}\right.
\end{eqnarray}
where $c'_i$ is a constant determined by $A, \wt{\gamma}$ and $\ell$, and independent of $r$.

\mslit
When $i=0$, by Lemma~\ref{lem: N-l}(3), we have   
\begin{eqnarray}
N^{pr}_{\ell+r}(\pi^\gamma, \pi^\gamma) &=& N^{pr}_{\ell+r}(\pi^v, \pi^v) = q^{rn(2n+1)}N_\ell^{pr}(\pi^v, \pi^v) \nonumber\\
&=& d_0\, q^{rn(2n+1)}, \label{N-ell+r0}
\end{eqnarray}
where $d_0 = N_\ell^{pr}(\pi^v, \pi^v)$ is a constant independent of $r$.
When $1 \leq i \leq n-1$, by Lemma~\ref{lem: N-l}, we obtain
\begin{eqnarray}
N^{pr}_{\ell+r}(\pi^\gamma,\pi^\gamma) &=& 
q^{(4(\ell+r)+2)i(n-i)+2i\abs{v}}\cdot q^{4\ell i^2}N_r^{pr}(\pi^\rho, \pi^\rho) \cdot q^{r(n-i)(2n-2i+1)}N_\ell^{pr}(\pi^v,\pi^v) \nonumber\\
&=& \label{N-(l+r)gamma}
q^{2i\abs{v} + 2i(n-i)+4\ell in +r(n-i)(2n+2i+1)} N_\ell^{pr}(\pi^v, \pi^v)N_r^{pr}(\pi^\rho,\pi^\rho)\nonumber\\
&=& \label{N-ell+r1}
d_i \, q^{r(n-i)(2n+2i+1)}N_r^{pr}(\pi^\rho,\pi^\rho),
\end{eqnarray}
where $d_i$ is a constant determined by $\wt{\gamma}$, and independent of $r$.

\mmslit
When $i=n$, 
\begin{eqnarray} \label{i=n}
N^{pr}_{\ell+r}(\pi^\gamma, \pi^\gamma) = N^{pr}_{\ell+r}(\pi^\ell\pi^\rho, \pi^\ell\pi^\rho) = q^{4\ell n^2}N_r^{pr}(\pi^\rho, \pi^\rho).
\end{eqnarray}
Now the result follows from \eqref{calS big-r}, \eqref{N-ell+r0}, \eqref{N-ell+r1} and \eqref{i=n}.
\qed

\slit
For formal power series $R_1(X)$ and $R_2(X)$ of $X$, we write $R_1(X) \sim R_2(X)$ if $R_1(X)-R_2(X)$ is a polynomial of $X$. By definition \eqref{P-wt-gamma} and Lemma~\ref{lem: ratio S-N}, we have the following.

\begin{prop} \label{prop: P-wt-gam}
Assume that $\wt{\gamma} \in \Lam_{n,2\ell}^+$ and set $i = m_{\wt{\gamma}}(2\ell)$. Then it holds
\begin{eqnarray}
P_{\wt{\gamma}}(A; X) &\sim&
\left\{\begin{array}{ll}  
a_0 \cdot \sum_{r=0}^\infty \left(q^{n(2n-2m-1)}X\right)^r & \mbox{if }i = 0\\
a_i \cdot \sum_{r=0}^\infty \left(q^{n(2n-2m-1)+i(2i-2m+1)}X\right)^r \displaystyle{\sum_{\rho \in \Lam_{i,2r}^+}}\, \dfrac{\calS_r(A, \pi^\rho)}{N_r^{pr}(\pi^\rho,\pi^\rho)} & \mbox{if }i \geq 1,
\end{array}\right. \nonumber
\end{eqnarray}
where $a_i$ is  the same constant as in Lemma~\ref{lem: ratio S-N}. \\ 
Especially, if $i=0$, then $P_{\wt{\gamma}}(A; X) \cdot (1-q^{n(2n-2m-1)}X)$ is a polynomial of $X$.
\end{prop}

\slit
We consider the case $\wt{\gamma} \in \Lam_{n, 2\ell}$ with $i = m_{\wt{\gamma}}(2\ell) \geq 1$. By Proposition~\ref{prop: P-wt-gam}, we have to calculate the sum for $\rho \in \Lam_{i, 2r}^+$.
When $r = 0$, we have
\begin{eqnarray*}
\Lam(\wt{\lam}, 0) = \{\wt{\lam}\}, \quad  \Lam_{i, 0}^+ = \{ \bf0\},
\end{eqnarray*}
which gives the constant term of $P_{\wt{\gamma}}(A; X)$. Hence it suffices to consider the sum for $r \geq 1$ in \eqref{P-wt-gamma}.  
Hereafter we assume that $r \geq 1$ and we decompose $\Lam_{i, 2r}^+$ as follows. 

\begin{eqnarray}
&&
\Lam_{i, 2r}^+ = \bigsqcup_{j=0}^{i-1}\, \Lam^{(j)}, \quad  \Lam^{(j)} = \Lam^{(j+)} \sqcup \Lam^{(j-)}, \quad \Lam^{(j-)} = \emptyset \; \mbox{ unless } i \equiv j \pmod{2}, \label{decomp Lam}\\
&&
\Lam^{(0+)} =\set{(2b)^i}{0 \leq b \leq r}, \nonumber\\
&&
\Lam^{(0-)} = 
\set{(2b+1)^i}{0 \leq b \leq r-1}  \mbox{ if } i \equiv 0\pmod{2},  \nonumber\\
&&
\Lam^{(j+)} = \set{(2b)^i+(\beta, 0^{i-j}) \in \Lam_{i,2r}^+}{0 \leq b \leq r-1, \; \beta \in \Lam_{j, 2(r-b)}^+, \; \beta_j \geq 1} \label{j-pm}\mbox{ if } j \geq 1,\nonumber\\
&&
\Lam^{(j-)} = 
\set{((2b)^j+\beta)(2b-1)^{i-j} \in \Lam_{i,2r}^+}{ 1 \leq b \leq r,\; \beta \in \Lam_{j,2(r-b)}^+}  \mbox{ if } j \geq 1, \, i \equiv j\pmod{2}. \nonumber 
\end{eqnarray}
We set for $1 \leq i \leq n$, according to the summation range 
\begin{eqnarray} 
Q_i (X) &:=& \label{Q i}
\sum_{r\geq 1} X^r \sum_{\rho \in \Lam^+_{i,2r}}\, \dfrac{S_r(A, \pi^\rho)}{N_r^{pr}(\pi^\rho, \pi^\rho)} = \sum_{j=0}^{i-1} Q_i^{(j)}(X), \nonumber\\ 
Q_i^{(j)}(X) &:=& \nonumber
\sum_{r\geq 1} X^r \sum_{\rho \in \Lam^{(j)}}\, \dfrac{S_r(A, \pi^\rho)}{N_r^{pr}(\pi^\rho, \pi^\rho)} = \sum_{k = +, -}\, Q_i^{(jk)}(X), \nonumber\\
Q_i^{(jk)}(X) &:=& \label{Q i-jk}
\sum_{r\geq 1} X^r \sum_{\rho \in \Lam^{(jk)}}\, \dfrac{S_r(A, \pi^\rho)}{N_r^{pr}(\pi^\rho, \pi^\rho)}, \quad(0 \leq j \leq i-1, \; k = \pm). \label{Q i-jk}
\end{eqnarray}  
Then, the assertion for the case $i = m_{\wt{\gamma}}(2\ell) \geq 1$ in Proposition~\ref{prop: P-wt-gam} can be written as follows. 
\begin{eqnarray} \label{P-g-Q}
P_{\wt{\gamma}}(A, X) \sim a_i Q_i(q^{n(2n-2m-1)+i(2i-2m+1)}X), 
\end{eqnarray}
where $a_i$ is the same constant as in Proposition~\ref{prop: P-wt-gam}.

\vspace{2cm}
\Section{Calculation of $Q_i(X)$ and the proof of Theorem~\ref{th: KitaokaS}}

{\bf 5.1.} Take $\rho = (2b+1)^i \in \Lam^{(0-)}, \; 0 \leq b \leq r-1$, for even $i$. Then $\pi^\rho = \pi^bH_*$, with $*= i/2$, and 
\begin{eqnarray*}
\frac{\calS_r(A, \pi^\rho)}{N_r^{pr}(\pi^\rho,\pi^\rho)}
&=& q^{4bi(m-i)} \frac{\calS_{r-b}(A, H_*)}{N_{r-b}^{pr}(H_*,H_*)}.
\end{eqnarray*}
Thus we obtain
\begin{eqnarray*}
Q_i^{(0-)}(X) &=&
\sum_{r \geq 1}\, X^r \sum_{b=0}^{r-1}\, q^{4bi(m-i)}\dfrac{\calS_{r-b}(A, H_*)}{N_{r-b}^{pr}(H_*, H_*)} = 
\sum_{b \geq 0}\, \sum_{r \geq b+1}\,  X^r q^{4bi(m-i)}\dfrac{\calS_{r-b}(A, H_*)}{N_{r-b}^{pr}(H_*,H_*)}\nonumber \\
&=&
\sum_{b \geq 0} \sum_{r \geq 1} X^{r+b}q^{4bi(m-i)}\dfrac{\calS_r(A, H_*)}{N_r^{pr}(H_*, H_*)} =
\sum_{b \geq 0} (q^{4i(m-i)}X)^b \sum_{r \geq 1} X^r\dfrac{\calS_r(A, H_*)}{N_r^{pr}(H_*, H_*)}.
\end{eqnarray*}
Together with Lemma~\ref{lem: low ratio S-N}, we see  
\begin{eqnarray}
(1-q^{4i(m-i)}X) \cdot Q_i^{(0-)}(X) &\sim& c\, \sum_{r \geq 1} (q^{i(2m-2i-1)}X)^r,
\end{eqnarray}
where $c$ is a constant independent of $r$, which equals to $c_1$ in Lemma~\ref{lem: low ratio S-N}. 
Thus we see that
\begin{eqnarray}   \label{Q-i-0-}
(1-q^{4i(m-i)}X)(1-q^{i(2m-2i-1)}X) \cdot Q_i^{(0-)}(X)\; \mbox{is a polynomial in $X$}
\end{eqnarray}

\slit
{\bf 5.2.} Take $\rho = (2b)^i \in \Lam^{(0+)}, \; 0 \leq b \leq r$. Then $\pi^\rho = \pi^b1_i$, and
\begin{eqnarray}
\frac{\calS_r(A, \pi^\rho)}{N_r^{pr}(\pi^\rho,\pi^\rho)} &=&  
\left\{\begin{array}{ll}
q^{4bi(m-i)} \frac{\calS_{r-b}(A, 1_i)}{N^{pr}_{r-b}(1_i, 1_i)} & \mbox{if } b < r,\\
q^{4ri(m-i)}\frac{1}{w_i(q^{-2})} & \mbox{if } b = r.
\end{array}\right. \nonumber\\
\end{eqnarray}
Thus we obtain
\begin{eqnarray}
Q_i^{(0+)}(X) &=&
\sum_{r \geq 1}\, X^r \left(\sum_{b=0}^{r-1}\, q^{4bi(m-i)}\dfrac{\calS_{r-b}(A, 1_i)}{N_{r-b}^{pr}(1_i,1_i)}+ c_1 q^{4ri(m-i)} \right)\nonumber \\
&=& 
\dfrac{c_1 X}{1-q^{4i(m-i)}X} + \sum_{b \geq 0}\, \sum_{r \geq b+1}\,  X^r q^{4bi(m-i)}\dfrac{\calS_{r-b}(A, 1_i)}{N_{r-b}^{pr}(1_i,1_i)}\nonumber \\
&=&
\dfrac{c_1 X}{1-q^{4i(m-i)}X} + \sum_{b \geq 0} \sum_{r \geq 1} X^{r+b}q^{4bi(m-i)}\dfrac{\calS_r(A, 1_i)}{N_r^{pr}(1_i,1_i)}\nonumber \\
&=&
\dfrac{c_1 X}{1-q^{4i(m-i)}X} +\sum_{b \geq 0} (q^{4i(m-i)}X)^r \sum_{r \geq 1} X^r\dfrac{\calS_r(A, 1_i)}{N_r^{pr}(1_i,1_i)},
\end{eqnarray}
where $c_1= w_i(q^{-2})^{-1}$ is a constants independent of $r$.
Together with Lemma~\ref{lem: low ratio S-N}, we see  
\begin{eqnarray}
(1-q^{4i(m-i)}X) \cdot Q_i^{(0+)}(X) &\sim& c \sum_{r \geq 1} (q^{i(2m-2i-1)}X)^r,
\end{eqnarray}
where $c$ is a constant independent of $r$, which equals to $c_0$ in Lemma~\ref{lem: low ratio S-N}. Thus we see that
\begin{eqnarray} \label{Q-i-0+}
(1-q^{4i(m-i)}X)(1-q^{i(2m-2i-1)}X) \cdot Q_i^{(0+)}(X)\; \mbox{is a polynomial in $X$},
\end{eqnarray}
where the result is the same for $Q_i^{(0-)}(X)$ (\eqref{Q-i-0-}). Since $Q_i^{(0)}(X) = Q_i^{(0+)}(X)+ Q_i^{(0-)}(X)$, we have
\begin{eqnarray} \label{Q-i j=0}
(1-q^{4i(m-i)}X)(1-q^{i(2m-2i-1)}X) \cdot Q_i^{(0)}(X)\; \mbox{is a polynomial in $X$}.
\end{eqnarray}

\mslit
Assume that $\wt{\gamma} \in \Lam_{n, 2\ell}$ and $m_{\wt{\gamma}}(2\ell) = 1$. 
Then $Q_1(X) = Q_1^{(0)}(X)$ and \\$P_{\wt{\gamma}}(X) \sim a_1 Q_1(q^{n(2n-2m-1)+(3-2m)}X)$  (cf. \eqref{Q i-jk} and \eqref{P-g-Q}). Thus we obtain the following. 

\mslit
\begin{prop}\label{prop: i=1}
Assume $\wt{\gamma} \in \Lam_{n, 2\ell}$ and $m_{\wt{\gamma}}(2\ell) = 1$. Then one has
\begin{eqnarray*}
&&
(1-q^{4(m-1)}X)(1-q^{2m-3}X) \cdot Q_1(X) \mbox{ is a polynomial in $X$},\\
&&
(1-q^{n(2n-2m-1)+2m-1}X)(1-q^{n(2n-2m-1)}X) \cdot P_{\wt{\gamma}}(A, X) \mbox{ is a polynomial in $X$}.
\end{eqnarray*}
\end{prop}

\mslit
{\bf 5.3.} Assume that $j \geq 1$ and take $\rho \in \Lam^{(j-)}(\subset \Lam_{i,2r})$. Then we may write (cf. \eqref{decomp Lam})
\begin{eqnarray*}\pi^\rho = \pi^b\pi^\beta \bot \pi^{b-1}H_*, \quad (b \geq 1, \; \beta \in \Lam_{j,2(r-b)}^+, \; *=(i-j)/2).
\end{eqnarray*}  
Then we have
\begin{eqnarray}
\calS_r(A, \pi^\rho)&=&
\left\{\begin{array}{ll}
q^{4bjm+4(b-1)(i-j)m}\calS_{r-b}(A, \pi^\beta)\calS_{r-b+1}(A, H_*) & \mbox{if } b < r\\[2mm]
q^{4rjm + 4(r-1)(i-j)m} \calS_1(A, H_*) & \mbox{if } b = r, \end{array}\right.\label{decomp S, N}\\
N_r^{pr}(\pi^\rho,\pi^\rho)
&=&
\left\{\begin{array}{ll}
q^{(4r+2)j(i-j)+2j(2b-1)(i-j)+4bj^2+4(b-1)(i-j)^2 }N^{pr}_{r-b}(\pi^\beta, \pi^\beta) N_{r-b+1}^{pr}(H_*,H_*) & \mbox{if } b < r\\
q^{(4r+2)j(i-j)+2j(2r-1)(i-j)}q^{4rj^2+4(r-1)(i-j)^2}w_j(q^{-2})N_1^{pr}(H_*,H_*) & \mbox{if } b=r. 
\end{array}\right. \nonumber
\end{eqnarray}
Then we obtain
\begin{eqnarray}
\dfrac{\calS_r(A, \pi^\rho)}{N_r^{pr}(\pi^\rho,\pi^\rho)}
&=&
\left\{\begin{array}{ll}
q^{rf+bg+h}\dfrac{\calS_{r-b}(A, \pi^\beta)}{N_{r-b}^{pr}(\pi^\beta,\pi^\beta)} \dfrac{\calS_{r-b+1}(A, H_*)}{N_{r-b+1}^{pr}(H_*, H_*)} & \mbox{if } b < r\\[3mm]
c_1 q^{4ri(m-i)}\dfrac{\calS_1(A, H_*)}{N_1^{pr}(H_*, H_*)} & \mbox{if } b = r
\end{array}\right.
\label{with f,g,h}
\end{eqnarray}
where $f, g, h$ are polynomials of $i, j, m$, rewrited from \eqref{decomp S, N} and independent of $r$ nor $b$, and $c_1$ is a constant independent of $r$. Thus we obtain

\begin{eqnarray}
\lefteqn{Q_i^{(j-)}(X)}\\
 &=&
\sum_{r \geq 1}X^r \left(\sum_{b=1}^{r-1} \sum_{\beta\in \Lam_{j,2(r-b)}^+}\, q^{rf+bg+h}\dfrac{\calS_{r-b}(A, \pi^\beta)}{N_{r-b}^{pr}(\pi^\beta,\pi^\beta)} \dfrac{\calS_{r-b+1}(A, H_*)}{N_{r-b+1}^{pr}(H_*, H_*)}+c_1q^{4ri(m-i)}\dfrac{\calS_1(A, H_*)}{N_1^{pr}(H_*, H_*)}\right)\nonumber \\
&=&
q^h\sum_{r\geq 1}\, \sum_{b=1}^{r-1}\, X^r\sum_{\beta\in \Lam_{j,2(r-b)}^+}\, q^{rf+bg}\dfrac{\calS_{r-b}(A, \pi^\beta)}{N_{r-b}^{pr}(\pi^\beta,\pi^\beta)} \dfrac{\calS_{r-b+1}(A, H_*)}{N_{r-b+1}^{pr}(H_*, H_*)} + c_1 \sum_{r \geq 1}\, (q^{4i(m-i)}X)^r\dfrac{\calS_1(A, H_*)}{N_1^{pr}(H_*, H_*)}\nonumber\\
&=&
q^h\sum_{b \geq 1}\, \sum_{r \geq b+1}\, X^r\sum_{\beta\in \Lam_{j,2(r-b)}^+}\, q^{rf+bg}\dfrac{\calS_{r-b}(A, \pi^\beta)}{N_{r-b}^{pr}(\pi^\beta,\pi^\beta)} \dfrac{\calS_{r-b+1}(A, H_*)}{N_{r-b+1}^{pr}(H_*, H_*)} + \dfrac{c_2 X}{1-q^{4i(m-i)}X}\nonumber\\
&=&
q^h\sum_{b \geq 1}\, \sum_{r \geq 1}\, X^{r+b}\sum_{\beta\in \Lam_{j,2r}^+}\, q^{(r+b)f+bg}\dfrac{\calS_{r}(A, \pi^\beta)}{N_r^{pr}(\pi^\beta,\pi^\beta)} \dfrac{\calS_{r+1}(A, H_*)}{N_{r+1}^{pr}(H_*, H_*)} + \dfrac{c_2 X}{1-q^{4i(m-i)}X)}\nonumber\\
&=&
q^h \sum_{b \geq 1}\, (q^{f+g}X)^b \underbrace{\sum_{r \geq 1}\, (q^f X)^r 
\sum_{\beta\in \Lam_{j,2r}^+}\, \dfrac{\calS_r(A, \pi^\beta)}{N_r^{pr}(\pi^\beta,\pi^\beta)} \dfrac{\calS_{r+1}(A, H_*)}{N_{r+1}^{pr}(H_*, H_*)}}_{(\sharp)} + \dfrac{c_2 X}{1-q^{4i(m-i)}X},\nonumber
\end{eqnarray}
where $c_2$ is a constant independent of $r$ nor $b$.
By Lemma~\ref{lem: low ratio S-N}, we see
\begin{eqnarray}
(\sharp) &\sim& c_3 \sum_{r\geq 1} (q^{f+(i-j)(2m-2i+2j-1)}X)^r \sum_{\beta\in \Lam_{j,2r}^+}\,\dfrac{\calS_r(A, \pi^\beta)}{N_r^{pr}(\pi^\beta,\pi^\beta)}\nonumber \\
&=&
c_3 \, Q_j(q^{(i-j)(2m-2i-2j-1)}X),
\end{eqnarray}
where, by \eqref{decomp S, N}, 
\begin{eqnarray*}
f=-4j(i-j), \quad f + g = -4j(i-j) + 4(m-i^2+ij-j^2) = 4i(m-i). 
\end{eqnarray*}
Thus, together with  a suitable constant $c$, we obtain

\begin{eqnarray} \label{Q j-}
(1-q^{4i(m-i)}X) \cdot Q_i^{(j-)}(X) \sim c\, X \cdot Q_j(q^{(i-j)(2m-2i-2j-1)}X) \quad (j \geq 1).
\end{eqnarray}

\mslit
{\bf 5.4.} 
Assume that $j \geq 1$ and we decompose $\Lam^{(j+)} (\subset \Lam_{i, 2r}^+)$ as follows ($r \geq 1$ is assumed):
\begin{eqnarray}
&&
\Lam^{(j+)} = \bigsqcup_{k=0}^{[j/2]} \Lam^{(j+,k)} \nonumber \\
&&
\mbox{for $0 \leq k \leq [j/2]$ and $j-2k >0$, } \nonumber\\
&& \qquad
\Lam^{(j+,k)} = \set{((2b+2)^{j-2k}+\gamma)(2b+1)^{2k}(2b)^{i-j} \in \Lam_{i,2r}}{0\leq b \leq r-1, \; \gamma \in \Lam_{j-2k, 2(r-b-1)}^+}, \nonumber\\
&&
\mbox{for even $j$ and $k = j/2$, }\nonumber\\
&& \qquad
\Lam^{(j+, j/2)} = \set{(2b+1)^{j}(2b)^{i-j} \in \Lam_{i,2r}}{0\leq b \leq r-1},\label{Lj+decomp}
\end{eqnarray}
and define $Q_i^{(j+,k)}(X)$ as before. Assume that $k \geq 1, \, j-2k>0$ and take $\rho \in \Lam^{j+,k}$ and write it as in \eqref{Lj+decomp}. 
Then $\pi^\rho = \pi^{b+1}\pi^\gamma \bot \pi^b H_k \bot \pi^b 1_{i-j}$, and we have, by Lemma~\ref{prop: calS(A,C)} and Lemma~\ref{lem: N-l}
\begin{eqnarray}
\calS_r(A, \pi^\rho)&=&
q^{4bmi}q^{4m(j-2k)} \calS_{r-b-1}(A, \pi^\gamma) \calS_{r-b}(A, H_k) \calS_{r-b}(A, 1_{i-j}) ,\nonumber\\
N_r^{pr}(\pi^\rho,\pi^\rho)&=&
q^{4bi^2} q^{(4(r-b)+2)j(i-j) }q^{(4(r-b)+2)2k(j-2k)+4k(j-2k)}\nonumber\\
&&
\times N_{r-b}^{pr}(\pi\pi^\gamma, \pi\pi^\gamma)N_{r-b}^{pr}(H_k, H_k) N_{r-b}^{pr}(1_{i-j},1_{i-j}). \nonumber
\end{eqnarray}
If $b = r-1$, we see  
\begin{eqnarray}
\dfrac{S_r(A, \pi^\rho)}{N_r^{pr}(\pi^\rho, \pi^\rho)} &=& c_1 \cdot q^{4ri(m-i)}\label{S-N j+k, 0},
\end{eqnarray}
where $c_1$ is a constant independent of $r$. While, if $0 \leq b \leq r-2$, we see
\begin{eqnarray}
\dfrac{S_r(A, \pi^\rho)}{N_r^{pr}(\pi^\rho, \pi^\rho)} &=& 
c_2 \cdot q^{rf+bg} \dfrac{\calS_{r-b-1}(A, \pi^\gamma)\calS_{r-b}(A, H_k) \calS_{r-b}(A, 1_{i-j})}{N_{r-b-1}^{pr}(\pi^\gamma,\pi^\gamma)N_{r-b}^{pr}(H_k, H_k) N_{r-b}^{pr}(1_{i-j},1_{i-j})}, \nonumber \\
\label{S-N j+k, main}
\end{eqnarray}
where
$f= -4(j(i-j) +2k(j-2k)), \quad g = 4(mi-i^2+j(i-j)+2k(j-2k))$,  and $c_2$ is a constant independent of $r$ nor $b$. In this case $f+g=4i(m-i)$ also (cf. the end of {\bf 5.3}). 

\mslit
Now we have
\begin{eqnarray}
\lefteqn{Q_i^{(j+,k)}(X) = \sum_{r\geq 1}\, X^r \sum_{\rho \in \Lam^{(j+,k)}}\, \dfrac{S_r(A, \pi^\rho)}{N_r^{pr}(\pi^\rho, \pi^\rho)}} \nonumber\\
&=&
\sum_{r\geq 1}X^r \left(\sum_{b=0}^{r-2}\left(\sum_{\gamma \in \Lam_{j-2k,2(r-b-1)}} c_2 q^{rf+bg} \dfrac{\calS_{r-b-1}(A, \pi^\gamma)\calS_{r-b}(A, H_k) \calS_{r-b}(A, 1_{i-j})}{N_{r-b-1}^{pr}(\pi^\gamma,\pi^\gamma)N_{r-b}^{pr}(H_k, H_k) N_{r-b}^{pr}(1_{i-j},1_{i-j})}\right)\right.\nonumber\\
&&
\qquad  + c_1q^{4ri(m-i)}\Big{)} \nonumber\\
&=&
c_2\sum_{b \geq 0}\sum_{r \geq b+2} X^rq^{rf+bg}\sum_{\gamma \in \Lam_{j-2k,2(r-b-1)}}\dfrac{\calS_{r-b-1}(A, \pi^\gamma)\calS_{r-b}(A, H_k) \calS_{r-b}(A, 1_{i-j})}{N_{r-b-1}^{pr}(\pi^\gamma,\pi^\gamma)N_{r-b}^{pr}(H_k, H_k) N_{r-b}^{pr}(1_{i-j},1_{i-j})} \nonumber\\
&&
\qquad + c_1 \sum_{r \geq 1}\, (q^{4i(m-i)}X)^r \nonumber\\
&=&
c_2 \sum_{b \geq 0}\sum_{r \geq 1}X^{r+b+1}q^{(r+b+1)f+bg}\sum_{\gamma \in \Lam_{j-2k,2r}}
\dfrac{\calS_r(A, \pi^\gamma)\calS_{r+1}(A, H_k) \calS_{r+1}(A, 1_{i-j})}{N_r^{pr}(\pi^\gamma,\pi^\gamma)N_{r+1}^{pr}(H_k, H_k) N_{r+1}^{pr}(1_{i-j},1_{i-j})} \nonumber\\
&&
\qquad + \frac{c_1 X}{1-q^{4i(m-i)}X}\nonumber\\
&=&
c_2 q^f X \sum_{b \geq 0} (q^{f+g}X)^b \underbrace{\sum_{r \geq 1}(q^fX)^r \sum_{\gamma \in \Lam_{j-2k,2r}}
\dfrac{\calS_r(A, \pi^\gamma)\calS_{r+1}(A, H_k) \calS_{r+1}(A, 1_{i-j})}{N_r^{pr}(\pi^\gamma,\pi^\gamma)N_{r+1}^{pr}(H_k, H_k) N_{r+1}^{pr}(1_{i-j},1_{i-j})}}_{(\sharp2)} \nonumber\\
&&
\qquad + \frac{c_1 X}{1-q^{4i(m-i)}X}. \nonumber
\end{eqnarray}
By Lemma~\ref{lem: low ratio S-N}, we see
\begin{eqnarray}
(\sharp2) &\sim & c_3 \sum_{r \geq 1}(q^{f+2k(2m-4k-1)+(i-j)(2m-2i+2j-1)}X)^r \sum_{\gamma \in \Lam_{j-2k,2r}} \dfrac{\calS_r(A, \pi^\gamma)}{N_r^{pr}(\pi^\gamma,\pi^\gamma)}\nonumber\\
&=&
c_3\,  Q_{j-2k}(q^{(i-j)(2m-2i-2j-1)+2k(2m-4j+4k-1)}X) \nonumber\\
&=&
c_3\,  Q_{j-2k}(q^{(i-j+2k)(2m-2i-2j+4k-1)}X).
\end{eqnarray}
Thus, for $1 \leq k \leq [\frac{j}{2}],\; j-2k>0$, together with a suitable constant $a_k$, which is new and independent of $r$, we obtain, , 
\begin{eqnarray} \label{Q j+k}
(1-q^{4i(m-i)}X)\cdot Q_i^{(j+,k)}(X) \sim a_k\, X\cdot Q_{j-2k}(q^{(i-j+2k)(2m-2i-2j+4k-1)}X).
\end{eqnarray}
 
\slit
Next we consider the case $k = 0$ and take $\rho \in \Lam^{(j+,0)}$. By \eqref{Lj+decomp}, $\pi^\rho = \pi^{b+1}\pi^\gamma \bot \pi^b 1_{i-j}, \gamma \in \Lam_{j,2(r-b-1)}^+,\; 0 \leq b \leq r-1$. 
When $b = r-1$, similar to the case $k \geq 1$, we see
\begin{eqnarray}
\frac{\calS_r(A, \pi^\rho)}{N_r^{pr}(\pi^\rho,\pi^\rho)}&=& c_4\cdot q^{4ri(m-i)},
\end{eqnarray}
where $c_4$ is a constant independent on $r$. 
While, if $0 \leq b \leq r-2$, 
\begin{eqnarray}
\frac{\calS_r(A, \pi^\rho)}{N_r^{pr}(\pi^\rho,\pi^\rho)}&=& 
\frac{q^{4bmi+4mj}\calS_{r-b-1}(A, \pi^\gamma)\calS_{r-b}(A, 1_{i-j})}
{q^{4bi^2+(4(r-b)+2)j(i-j) + 4j^2}N_{r-b-1}^{pr}(\pi^\gamma, \pi^\gamma)N_{r-b}^{pr}(1_{i-j},1_{i-j}) } \nonumber \\
&=&
c_5 \cdot q^{rf+bg} \, \dfrac{\calS_{r-b-1}(A, \pi^\gamma)}{N_{r-b-1}^{pr}(\pi^\gamma,\pi^\gamma)}
\dfrac{\calS_{r-b}(A, 1_{i-j})}{N_{r-b}^{pr}(1_{i-j}, 1_{i^j})},
\end{eqnarray}
where $f = -4j(i-j), \; g = 4(mi-i^2+j(i-j))$, and $c_5$ is a constant independent on $r$ nor $b$.
($f$ and $g$ are the same as in the case $k \geq 1$, if putting $k=0$.)
Thus we have
\begin{eqnarray}
\lefteqn{Q_i^{(j+, 0)}(X)} \nonumber\\
&=&
\sum_{r\geq 1}X^r \left(\sum_{b=0}^{r-2}\left(\sum_{\gamma \in \Lam_{j,2(r-b-1)}} c_5 q^{rf+bg} \dfrac{\calS_{r-b-1}(A, \pi^\gamma) \calS_{r-b}(A, 1_{i-j})}{N_{r-b-1}^{pr}(\pi^\gamma,\pi^\gamma) N_{r-b}^{pr}(1_{i-j},1_{i-j})}\right) + c_4q^{4ri(m-i)}\right) \nonumber\\
&=&
c_5\sum_{b \geq 0}\sum_{r \geq b+2} X^rq^{rf+bg}\sum_{\gamma \in \Lam_{j,2(r-b-1)}}\dfrac{\calS_{r-b-1}(A, \pi^\gamma) \calS_{r-b}(A, 1_{i-j})}{N_{r-b-1}^{pr}(\pi^\gamma,\pi^\gamma) N_{r-b}^{pr}(1_{i-j},1_{i-j})} + c_4 \sum_{r \geq 1}\, (q^{i(m-i)}X)^r\nonumber\\
&=&
c_5\sum_{b \geq 0}\, \sum_{r \geq 1} X^{r+b+1}q^{(r+b+1)f+bg}\sum_{\gamma \in \Lam_{j,2r}}\, \dfrac{\calS_r(A, \pi^\gamma) \calS_{r+1}(A, 1_{i-j})}{N_r^{pr}(\pi^\gamma,\pi^\gamma) N_{r+1}^{pr}(1_{i-j},1_{i-j})} + \frac{c_4 X}{1-q^{4i(m-i)}X}\nonumber\\
&=&
c_5q^fX \sum_{b\geq 0}\, (q^{f+g}X)^b \underbrace{\sum_{r \geq 1}\, (q^{f}X)^r  \sum_{\gamma \in \Lam_{j,2r}}\, \dfrac{\calS_r(A, \pi^\gamma) \calS_{r+1}(A, 1_{i-j})}{N_r^{pr}(\pi^\gamma,\pi^\gamma) N_{r+1}^{pr}(1_{i-j},1_{i-j})}}_{(\sharp3)} + \frac{c_4 X}{1-q^{4i(m-i)}X}, \nonumber\\
\end{eqnarray}
where $(\sharp3)$ is the same as $(\sharp2)$ by putting $k=0$. Hence 
\begin{eqnarray}
(\sharp3) &\sim & c_6 \sum_{r \geq 1}(q^{f+(i-j)(2m-2i+2j-1)}X)^r \sum_{\gamma \in \Lam_{j,2r}} \dfrac{\calS_r(A, \pi^\gamma)}{N_r^{pr}(\pi^\gamma,\pi^\gamma)}\nonumber\\
&=&
c_6\,  Q_{j}(q^{(i-j)(2m-2i-2j-1)}X).
\end{eqnarray}
Since $f+g = 4i(m-i)$, we obtain, with a suitable constant $a_0$, 
\begin{eqnarray} \label{Q j+0}
(1-q^{4i(m-i)}X)\cdot Q_i^{(j+,0)}(X) \sim a_0\, X\cdot Q_j(q^{(i-j)(2m-2i-2j-1)}X),
\end{eqnarray}
which is the same type as the case $Q_i^{(j-)}(X)$ (cf. \eqref{Q j-}) and \eqref{Q j+k} by putting $k=0$.

\slit
Finally we consider the case $j$ is even and $k= j/2$ and take $\rho \in \Lam^{(j+, \, j/2)}$. By 
\eqref{Lj+decomp}, $\pi^\rho = \pi^bH_{j/2} \bot \pi^b 1_{i-j}, \; 0 \leq b \leq r-1$.
When $b =r-1$, similar to the other cases, we see
\begin{eqnarray}
\frac{\calS_r(A, \pi^\rho)}{N_r^{pr}(\pi^\rho,\pi^\rho)}&=& c_7\cdot q^{4ri(m-i)},
\end{eqnarray}
where $c_7$ is a constant independent on $r$. 
While, if $0 \leq b \leq r-2$, 
\begin{eqnarray}
\frac{\calS_r(A, \pi^\rho)}{N_r^{pr}(\pi^\rho,\pi^\rho)}&=& 
\frac{q^{4bmi+4mj}\calS_{r-b}(A, H_{j/2})\calS_{r-b}(A, 1_{i-j})}
{q^{4bi^2+(4(r-b)+2)j(i-j) }N_{r-b}^{pr}(H_{j/2},H_{j/2}) N_{r-b}^{pr}(1_{i-j},1_{i-j}) } \nonumber \\
&=&
c_8 q^{rf+bg}\, \dfrac{\calS_{r-b}(A, H_{j/2})}{N_{r-b-1}^{pr}(H_{j/2},H_{j/2})}
\dfrac{\calS_{r-b}(A, 1_{i-j})}{N_{r-b}^{pr}(1_{i-j}, 1_{i^j})}.
\end{eqnarray}
where $f = -4j(i-j), \; g = 4(mi-i^2+j(i-j))$ and $c_8$ is a constant independent of $r$ nor $b$. 
Thus we have
\begin{eqnarray}
\lefteqn{Q_i^{(j+, \, j/2)}(X)} \nonumber\\
&=&
\sum_{r\geq 1}X^r \left(\sum_{b=0}^{r-2}\left( c_8 q^{rf+bg} \dfrac{\calS_{r-b}(A, H_{j/2}) \calS_{r-b}(A, 1_{i-j})}{N_{r-b}^{pr}(H_{j/2}. H_{j/2}) N_{r-b}^{pr}(1_{i-j},1_{i-j})}\right) + c_7q^{4ri(m-i)}\right) \nonumber\\
&=&
c_8\sum_{b \geq 0}\sum_{r \geq b+2} X^rq^{rf+bg}\dfrac{\calS_{r-b}(A, H_{j/2}) \calS_{r-b}(A, 1_{i-j})}{N_{r-b}^{pr}(H_{j/2}, H_{j/2}) N_{r-b}^{pr}(1_{i-j},1_{i-j})} + c_7 \sum_{r \geq 1}\, (q^{i(m-i)}X)^r\nonumber\\
&=&
c_8\sum_{b \geq 0}\, \sum_{r \geq 1} X^{r+b+1}q^{(r+b+1)f+bg}\, \dfrac{\calS_{r+1}(A, H_{j/2}) \calS_{r+1}(A, 1_{i-j})}{N_{r+1}^{pr}(H_{j/2}, H_{j/2}) N_{r+1}^{pr}(1_{i-j},1_{i-j})} + \frac{c_7 X}{1-q^{4i(m-i)}X}\nonumber\\
&=&
c_8q^fX \sum_{b\geq 0}\, (q^{f+g}X)^b \sum_{r \geq 1}\, c_9 (q^fX)^r q^{(j(m-2j-1)+(i-j)(m-2i+2j-1)r}+ \frac{c_7 X}{1-q^{4i(m-i)}X}\nonumber\\
&=&
cX\frac{1}{1-q^{4i(m-i)}X} \cdot \frac{X}{1-q^{i(2m-2i-1)}X}  + \frac{c_7 X}{1-q^{4i(m-i)}X},\nonumber
\end{eqnarray}
where $c, c_i$ are constants. Thus we see
\begin{eqnarray} \label{Qj+,j/2}
(1-q^{4i(m-i)}X)(1-q^{i(2m-2i-1)}X) \cdot Q_i^{(j+, \, j/2)}(X) \mbox{ is a polynomial},
\end{eqnarray}
which is the same type as the case $Q_i^{(0)}(X)$ (cf. \eqref{Q-i j=0}).

\slit
We recall that 
\begin{eqnarray*}
Q_i(X) = \sum_{j=0}^{i-1} Q_i^{(j)}(X) = Q_i^{(0)}(X)+ \sum_{j=1}^{i-1} \left(Q_i^{(j-)}(X)+\sum_{k=0}^{[j/2]}Q_i^{(j+,k)}(X)\right).
\end{eqnarray*}
Then, by \eqref{Q-i j=0},  \eqref{Q j-}, \eqref{Q j+k}, \eqref{Q j+0} and \eqref{Qj+,j/2}, 
\begin{eqnarray}
\lefteqn{(1-q^{4i(m-i)}X) (1-q^{i(2m-2i-1)}X)\cdot Q_i(X)}\nonumber\\
&\sim& 
(1-q^{4i(m-i)}X) (1-q^{i(2m-2i-1)}X) 
\left\{\sum_{1 \leq j \leq i-1}\, \left(Q_i^{(j-)}(X)+
\sum_{\scriptsize{\begin{array}{c} 0 \leq k \leq [j/2]\\ j-2k >0\end{array}}}\, Q_i^{(j+,k)}(X)\right)\right\} \nonumber\\
&\sim& 
(1-q^{i(2m-2i-1)}X) \left( \sum_{1 \leq j \leq i-1}\, \sum_{\scriptsize{\begin{array}{c} 0 \leq k \leq [j/2]\\ j-2k >0\end{array}}}\, c_{jk}\, X \cdot Q_{j-2k}(q^{(i-j+2k)(2m-2i-2j+4k-1)} X) \right) \nonumber\\
&\sim& 
(1-q^{i(2m-2i-1)}X) \left( \sum_{1 \leq j \leq i-1}\, c_j X \cdot Q_j(q^{(i-j)(2m-2i-2j-1)} X) \right), \label{Q-i}
\end{eqnarray}
where $c_{jk}$ and $c_j$ are constants. We adjust the variable by \eqref{P-g-Q}, i.e., we substitute $q^{n(2n-2m-1)+i(2i-2m+1)}X$ for $X$ in \eqref{Q-i}. Then we obtain, with suitable constants $c_j'$,
\begin{eqnarray}
\lefteqn{(1-q^{(n-i)(2n+2i-2m-1)}X) (1-q^{n(2n-2m-1)}X)\cdot Q_i(q^{n(2n-2m-1)+i(2i-2m+1)}X)}\nonumber\\ 
&\sim& 
(1-q^{n(2n-2m-1)}X) \left(\sum_{j=1}^{i-1}\, c'_{j}\, X \cdot Q_j(q^{n(2n-2m-1)+j(2j-2m+1)} X) \right) . \nonumber\\\label{adjust Q-i}
\end{eqnarray}
Thus, by induction, we obtain for $0 \leq i \leq n$, 
\begin{eqnarray*}
\prod_{j=0}^i\, (1-q^{(n-j)(2n+2j-2m-1)}X) \times Q_i(q^{n(2n-2m-1)+i(2i-2m+1)}X)\; \mbox{is a polynomial,}
\end{eqnarray*}
where, for $i = 0$, we set $Q_0(X) = \sum_{r \geq 1}\, X^r$ for convenience (cf. Proposition~\ref{prop: P-wt-gam}).   
Hence, by \eqref{P-g-Q}, we obtain for any $\wt{\gamma} \in \Lam_{n, 2\ell}^+$
\begin{eqnarray*}
\prod_{j=0}^n\, (1-q^{(n-j)(2n+2j-2m-1)}X) \times P_{\wt{\gamma}}(A, X) \; \mbox{is a polynomial.}
\end{eqnarray*}
Since $\Lam_{n,2\ell}^+$ is a finite set in Proposition~\ref{prop: P(B,A)},  we complete the proof of Theorem~\ref{th: KitaokaS}, i.e. we have proved that
\begin{eqnarray*} \label{result}
\prod_{i=0}^n\, (1-q^{(n-i)(2n+2i-2m-1)}X) \times P(B, A; X) \; \mbox{is a polynomial. \qquad \qed}
\end{eqnarray*}

\mslit
\begin{rem}
It is not easy to calculate local densities in general. By a direct calculation, we give a small example consistent with the assertion of Theorem~\ref{th: KitaokaS}.  
\begin{eqnarray*}
P(1, 1_2; X) = (1-q^{-2})\left(1 + \frac{1+q^{-1}}{1-q^{-3}}\cdot\frac{X}{1-X}-\frac{1+q^{-2}}{1-q^{-3}}\cdot\frac{q^{-4}X}{1-q^{-3}X}\right).
\end{eqnarray*}
\end{rem}
\vspace{2cm}
\setcounter{section}{5}
\Section{Appendix (Comparison chart)}
We note this appendix in order to compare the present main results on quaternion hermitian forms (Proposition~\ref{prop: F-tr of density}, Theorem~\ref{th: lin indep}, Theorem~\ref{th: KitaokaS}) with the previous results on symmetric forms and hermitian forms in \cite{BHS}, \cite{Kitaoka} and  \cite{Kitaoka2}.
Similar results for Kitaoka series of alternating forms are written in \cite{Kitaoka}, and it seems possible to prove similar results for the linear independence of local densities. 

\mmslit
The methods of the proofs are similar, but one has to verify case by case. The quaternion hermitian case is troublesome, since one has to calculate things on noncommutative algebra $D$ and the representatives of $K_n$-orbit decomposition are not so simple for calculation, even without unit part (cf. \eqref{pi-gamma}). Hence the simplest case should be the case of unramified hermitian forms, except for alternating forms.

\mslit
Let $k$ be a $\frp$-adic field, $\pi$ a prime element of $k$, and $q = \abs{\calo/\pi\calO_k}$. We assume that $q$\ is odd. 
Let $k'$ be a quadratic extension of $k$ with nontrivial $k$-isomorphism $*$ on $k'$, and we consider hermitian forms with respect to $*$. We have distinguished unramified hermitian forms (case (U)) and ramified hermitian forms (case (R)), and we write case (S) for symmetric forms. 
Set $\calO = \calO_k$ for case (S) and $\calO = \calO_{k'}$ for cases (U) and (R).
We define the spaces of symmetric, or hermitian forms (of size $n$) by
\begin{eqnarray}
&&
V_n = \left\{\begin{array}{ll}
\set{A \in M_n(k)}{{}^tA = A} & \mbox{for case (S)}\\
\set{A \in M_n(k')}{A^* = A} & \mbox{for cases (U) and (R)},
\end{array} \right. \nonumber \\
&&
X_n = \set{A \in V_n}{\mbox{non-degenerate}} \supset X_n^+ = \set{A\in X_n}{\mbox{integral}},
\end{eqnarray}
where $A^*$ is the conjugate transpose of $A$ with respect to $*$.
Local density $\mu(B, A)$ of $B \in X_n$ by $A \in X_m$ with $m \geq n$ is defined by 
\begin{eqnarray}
&&
\mu(B, A) = \lim_{\ell \rightarrow \infty}\, \frac{\sharp\set{\ol{v} \in M_{mn}(\calO/\pi^\ell\calO)}{A[v] -B \in M_n(\pi^\ell\calO)}}{q^{(*1)}},\\[2mm]
&&
\qquad 
\begin{array}{ll}
A[v] = {}^tv Av, \; (*1) = \ell n(m-(n+1)/2)  & \mbox{for case (S)}, \\[2mm]
A[v] = v^*Av, \; (*1) =  \ell n(2m-n) & \mbox{for cases (U) and (R)}.
\end{array} \nonumber
\end{eqnarray}

\mmslit
Let $\psi$ be an additive character on $k$ of conductor $\calo$, and define the pairing $\pair{\;}{\;}$ on $V_n\times V_n$ by $\pair{R}{T} = {\rm trace}(RT) \in k$ for $R, T \in V_n$.
For a function $f \in V_n$, its Fourier transform $\calF(f)$ is defined if the following integral is well-defined:
\begin{eqnarray}
\calF(f)(y) = \int_{V_n}\, f(x) \psi(-\pair{x}{y})dx.
\end{eqnarray}
For $S \in V_m$ and $T \in V_n$, we introduce the Gauss sum
\begin{eqnarray}
\calG(S, T) = \int_{M_{mn}(\calO)}\, \psi(\pair{S[v]}{T})dv.
\end{eqnarray}
The following lemma plays a key role (cf. Proposition~2.4, \cite[Lemma 1.2, Lemma 3.1]{BHS}). The convergent range for hermitian forms becomes better by using \cite[(13.4.2)]{Shimura}.

\begin{lem} 
For $A \in X_m$ and $B \in X_n$,  $\mu(B, A) = c_0 \calF(\calG(T,\; ))(B)$ if 
\begin{eqnarray*}
m \geq 2n+1 \mbox{ for case {\rm (S)}}, \quad m \geq 2n \mbox{ for case {\rm (U)}},
\quad m\geq 4n-1 \mbox{ for case {\rm (R)}},
\end{eqnarray*}
where $c_0 = 1$ for cases {\rm (S)} and {\rm (U)}, while $c_0= q^{-n(n-1)}$ for case {\rm (R)}.
\end{lem}

\mslit
We follow the notation (2.1) and (2.2), but we change the notation $h^r$ by
\begin{eqnarray}
h^r = \left\{\begin{array}{ll}
\twomatrix{0}{\pi^r}{\pi^r}{0} & \mbox{for cases (S) and (U)}\\
\twomatrix{0}{\varPi^r}{(-\varPi)^r}{0} & \mbox{for case (R)},
\end{array}\right.
\end{eqnarray}
where the prime element $\pi$ of $k$ is still a prime element of $k'$ for case (U), and $\varPi$ is a prime element of $k'$ for case (R) such that $\varPi^2 = \pi$. 
Then we have $H^\lam \in X_{2k}$ for $\lam \in \Gamma_k$.

\begin{thm} {\rm \cite[Theorem 2.4, Theorem 3.5]{BHS})}
Assume that $k, \ell, n \in \N$, and $r \in \Z, r \geq 0$ satisfy the condition $k \geq n$ and 
$m = 2k+r$ satisfies the condition of Lemma~6.1.

\noindent
Take an $S \in X_r^+$ if $r > 0$. Then one has the following.

\mmslit
{\rm (1)} For $\lam \in \Gamma_{k,\ell}$, $T \in X_n^+$, it satisfies
\begin{eqnarray*}
\mu(T, H^\lam \bot S) = \left\{\begin{array}{ll}
\sum_{\tau \in \Gamma_{n, \ell}}\, a_\tau q^{\pair{\what{\lam}}{\what{\tau}}} & \mbox{for cases {\rm (S)} and {\em (R)}}\\
\sum_{\tau \in \Gamma_{n, \ell}}\, a_\tau q^{2\pair{\what{\lam}}{\what{\tau}}} & \mbox{for case {\rm (U)}}, \end{array}\right. 
\end{eqnarray*}
where $a_\tau= a_\tau(k, S, T)$ is a rational constant independent of $\lam$.

\mmslit
{\rm (2)} As functions of $T \in X_n^+$, the set $\set{\mu(T, H^\lam \bot S)}{\lam \in \Gamma_{k,\ell}}$ spans a $\begin{pmatrix}n+\ell\\n \end{pmatrix}$-dimensional $\Q$-space, and the set $\set{\mu(T, H^\mu\bot H^{{\bf0}, k-n} \bot S)}{\mu \in \Gamma_{n,\ell}}$ forms a basis.
\end{thm}
 
\slit
Kitaoka series for $A \in X_m$ and $B \in X_n$ with $m \geq n$ is defined similarly in each case by
\begin{eqnarray}
P(B, A; X) = \sum_{r \geq 0}\, \mu(\pi^rB, A) X^r.
\end{eqnarray}

\begin{thm}{\rm ( Main theorems in \cite{Kitaoka}, \cite{Kitaoka2} )}
The power series $P(B, A; X)$ becomes a polynomial in $X$ if it is multiplied by the following polynomial for each case.

For case {\rm (S)} and even $m$, \; $\displaystyle{\prod_{i=0}^n\, (1 - (\eps_A q^{\frac12(n+i-m+1)})^{n-i}X)}$, 
where $\eps_A = \pm 1$ is given explicitly by $A$ in {\rm \cite{Kitaoka2}}. 

For case {\rm (S)} in general, \; $\displaystyle{(1 - X) \prod_{i=0}^{n-1}\, (1 - q^{(n-i)(n+i-m+1)}X^2)}.$

For case {\rm (U)},\; $\displaystyle{\prod_{i=0}^{n}\, (1 - (-1)^{m(n-i)}q^{(n-i)(n+i-m)}X)}$.

For case {\rm (R)} and $-1$ is square modulo $\frp$ or $m$ is even, \; $\displaystyle{\prod_{i=0}^{n}\, (1 - q^{(n-i)(n+i-m-1)}X)}$.

For case {\rm (R)} in general, \; $\displaystyle{(1 - X)\prod_{i=0}^{n-1}\, (1 - q^{2(n-i)(n+i-m-1)}X^2)}$.
\end{thm}
\vspace{2cm}
\noindent
{\it Acknowledgements}\\
This research is partially supported  by Grant-in-Aid for Scientific Research JP16K05081, and by RIMS, an International Joint Usage/Research Center located in Kyoto University.
The author is grateful to the referee for the careful reading and suitable comments.

\vspace{2cm}
\bibliographystyle{amsalpha}

\vspace{2cm}
\hfill \begin{minipage}{7cm}
Yumiko Hironaka\\
168-0072 Takaido-Higashi 3-9-10-730\\
Tokyo, JAPAN\\
e-mail: hironaka@waseda.jp
\end{minipage}%

\end{document}